\documentclass[twocolumn,oneside,reqno]{amsart}

\usepackage[textheight=250mm, textwidth=190mm]{geometry}
\usepackage{amssymb,amsmath,dsfont}
\usepackage{color,float}

\usepackage{hyperref}

\numberwithin{equation}{section}

\usepackage{extarrows}

\usepackage{enumitem}

\usepackage{caption}
\usepackage{array,booktabs}

\usepackage{tikz}
\usetikzlibrary{matrix,arrows}

\usepackage[norelsize,algo2e]{algorithm2e}

\def \mnc {\pgfmatrixnextcell}

\providecommand{\Natural}{\mathbb{N}}
\providecommand{\Integer}{\mathbb{Z}}
\providecommand{\Rational}{\mathbb{Q}}
\providecommand{\Real}{\mathbb{R}}

\providecommand{\Prob}{\mathbb{P}}
\providecommand{\Pois}{\operatorname{Pois}}
\providecommand{\algNorm}{\operatorname{N}}

\providecommand{\CyI}[1]{\Integer[\zeta_{#1}]}

\providecommand{\abs}[1]{\lvert#1\rvert}
\providecommand{\absAlt}[1]{\left\lvert#1\right\rvert}

\providecommand{\rprod}[2]{\Real^{#1} \times \Real^{#2}}

\providecommand{\intproj}{\pi_{\operatorname{int}}}

\definecolor{lgray}{gray}{0.35}

\providecommand{\genfigureobject}[3]{
\def\genfigPre{insertfigure:}
\def\genfigPost{#1}
\expandafter\def\csname\genfigPre\genfigPost\endcsname{%
  \begin{figure}[H]%
    \centering%
    #2%
    \caption{#3}%
    \label{fig:#1}%
  \end{figure}%
}}

\genfigureobject{z2lattice_vs_poisson}
  {\includegraphics[width=0.23\textwidth]{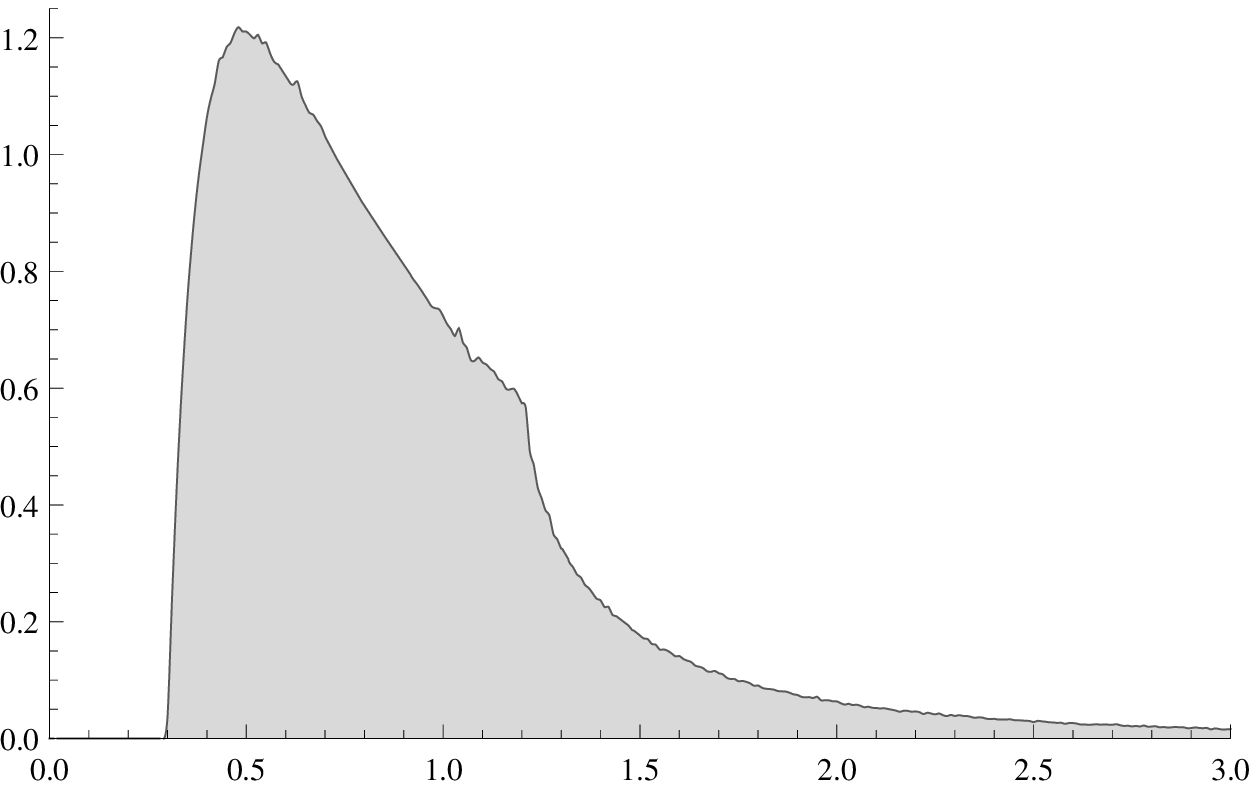} \hspace{1.0ex}
   \includegraphics[width=0.23\textwidth]{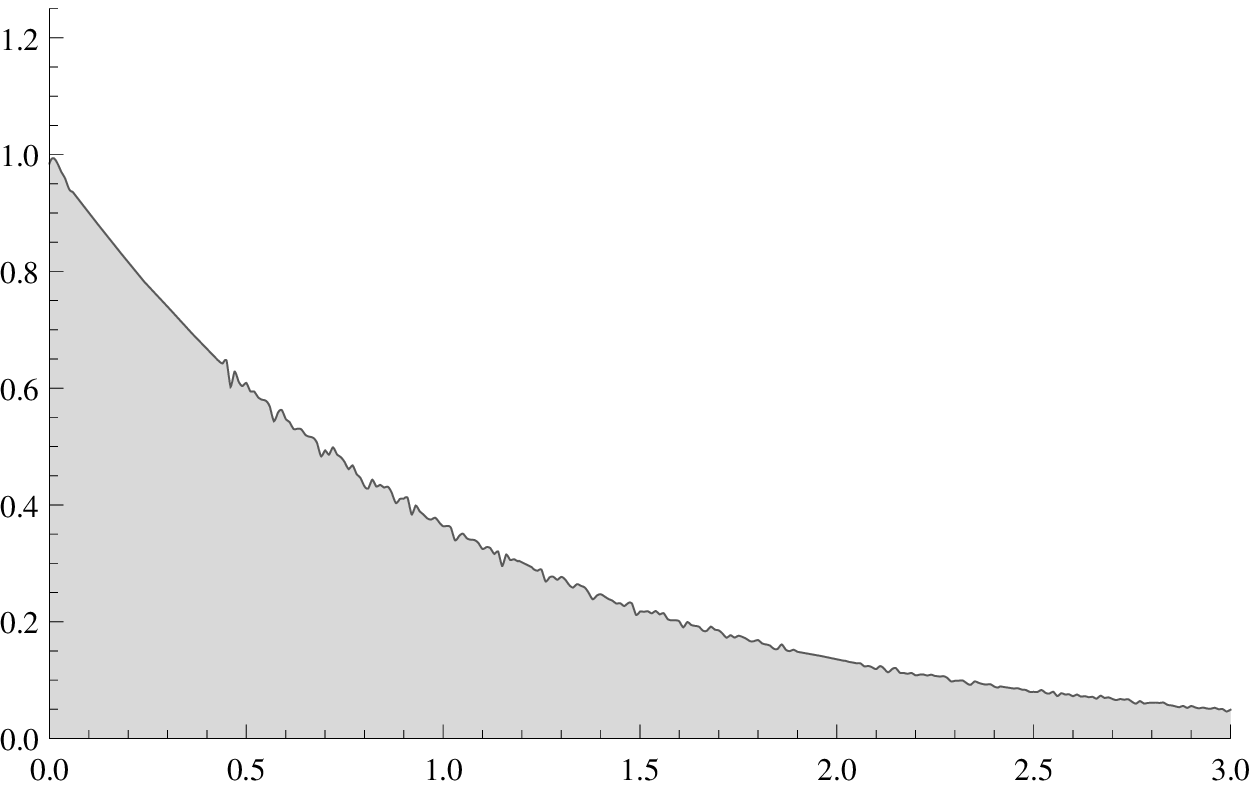}}
  {Asymptotic spacing distribution for $\Integer^2$ (left)
   and Poisson (right).}

\genfigureobject{ttt_inflrule_tile_a}
  {\includegraphics[width=0.092\textwidth]{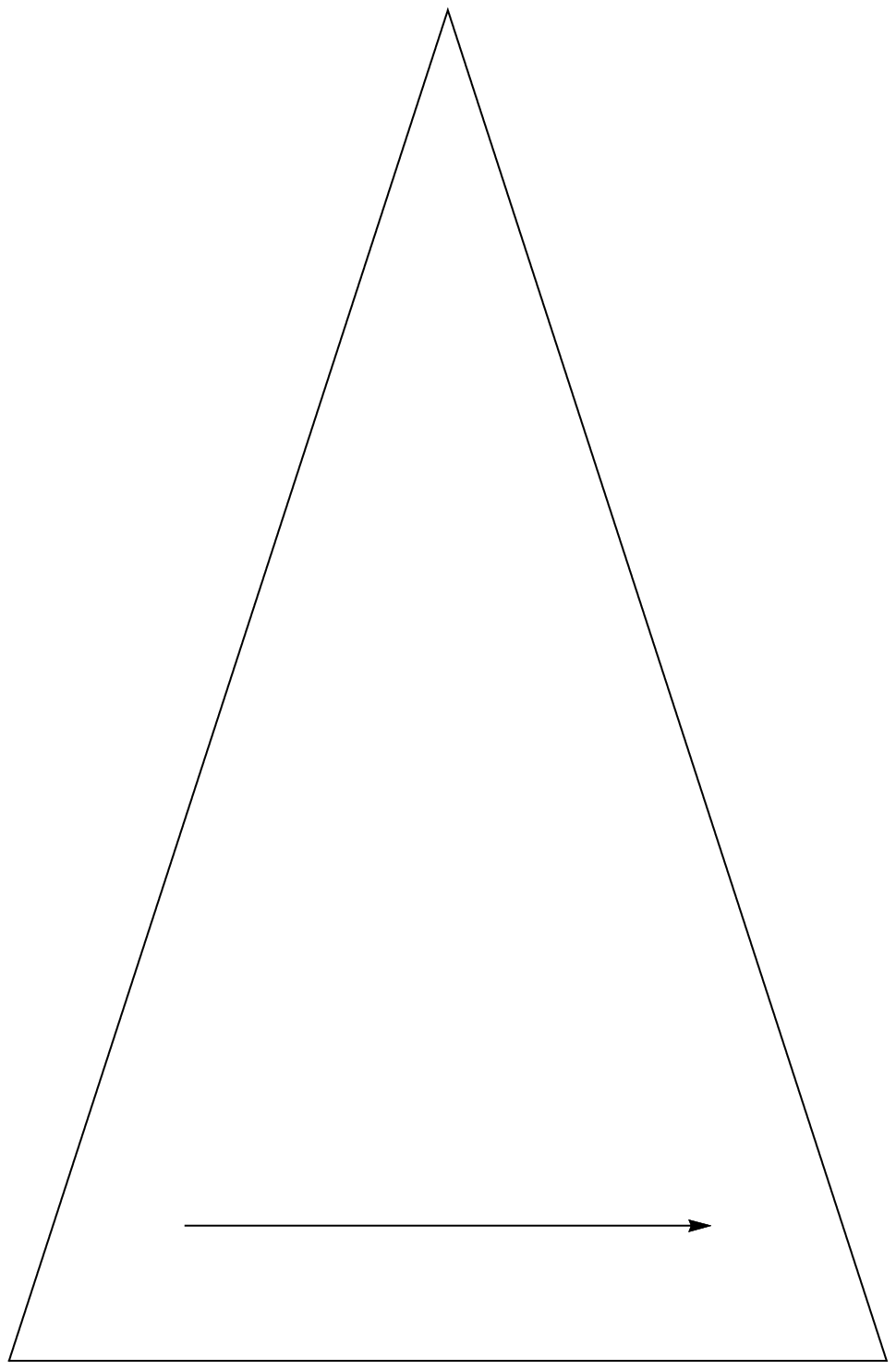} \hspace{1ex}
   \raisebox{0.06\textwidth}{$\xrightarrow{\hspace{9ex}}$} \hspace{1ex}
   \includegraphics[width=0.15\textwidth]{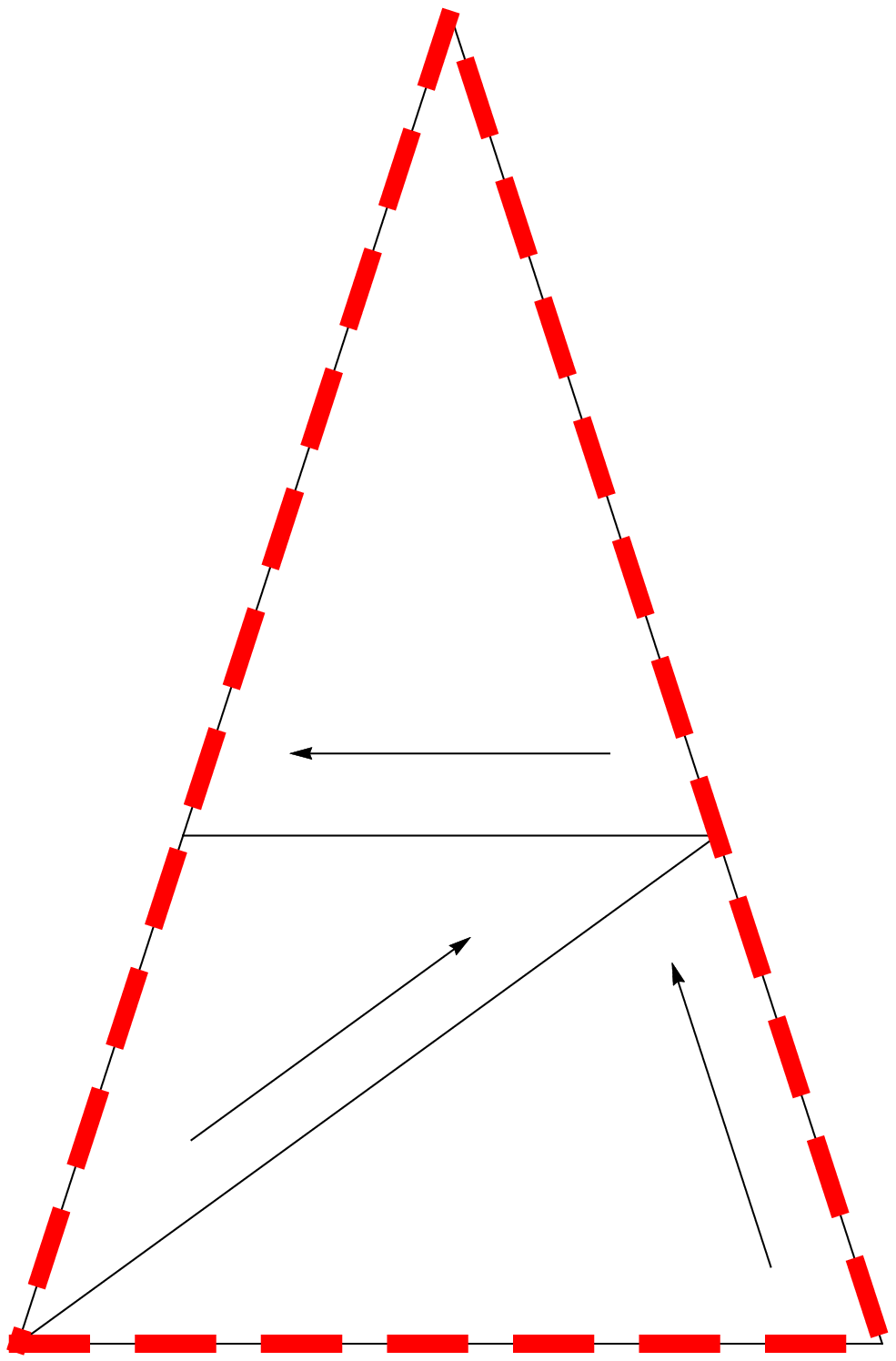}}
  {Tile A maps to 2$\times$A and 1$\times$B.}

\genfigureobject{ttt_inflrule_tile_b}
  {\includegraphics[width=0.092\textwidth]{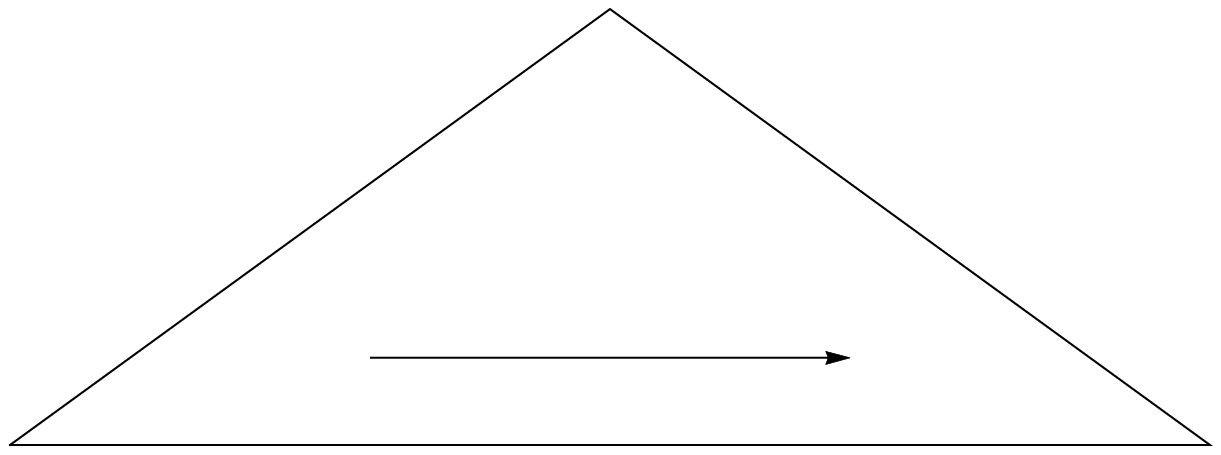} \hspace{1ex}
   \raisebox{0.015\textwidth}{$\xrightarrow{\hspace{9ex}}$} \hspace{1ex}
   \includegraphics[width=0.15\textwidth]{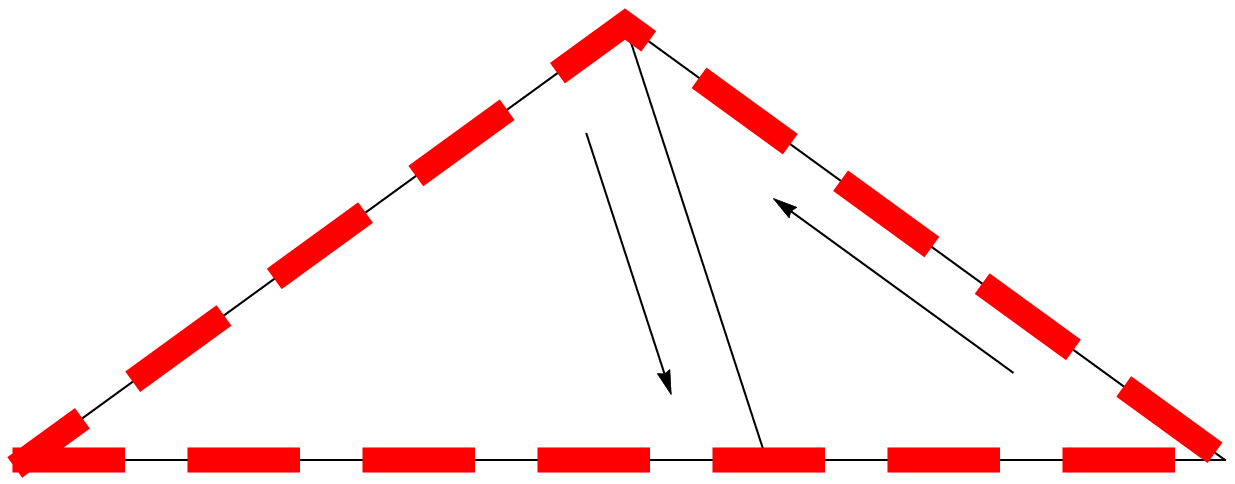}}
  {Tile B maps to 1$\times$A and 1$\times$B.}

\genfigureobject{ttt_apatch_infl5}
  {\includegraphics[width=0.36\textwidth]{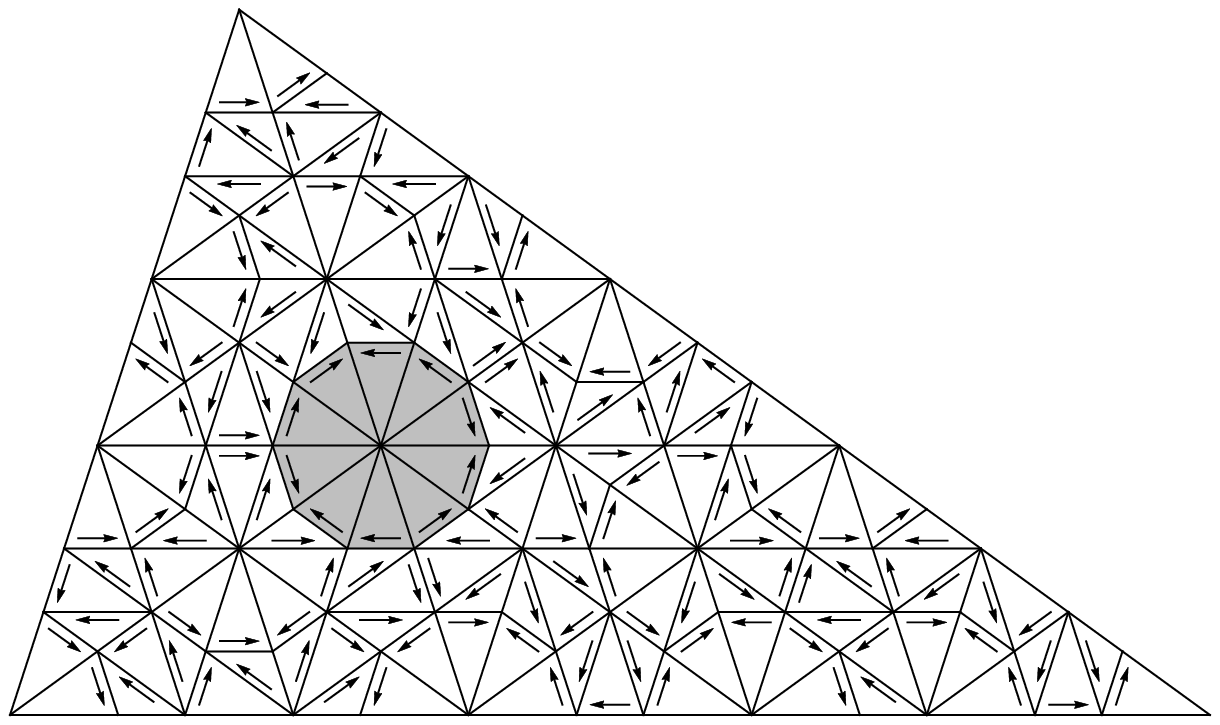}}
  {Patch produced from five inflation steps applied to the \textsf{TT} prototile \emph{A}.}

\genfigureobject{ab_inflrule_tile_a}
  {\includegraphics[width=0.066\textwidth]{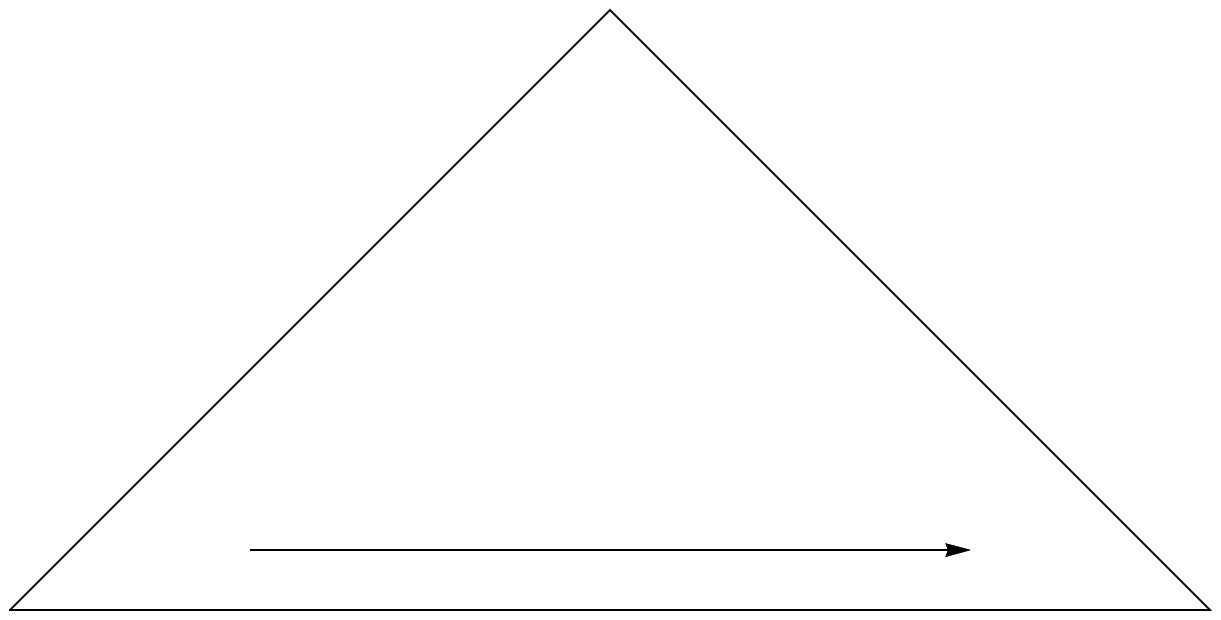} \hspace{1ex}
   \raisebox{0.015\textwidth}{$\xrightarrow{\hspace{9ex}}$} \hspace{1ex}
   \includegraphics[width=0.16\textwidth]{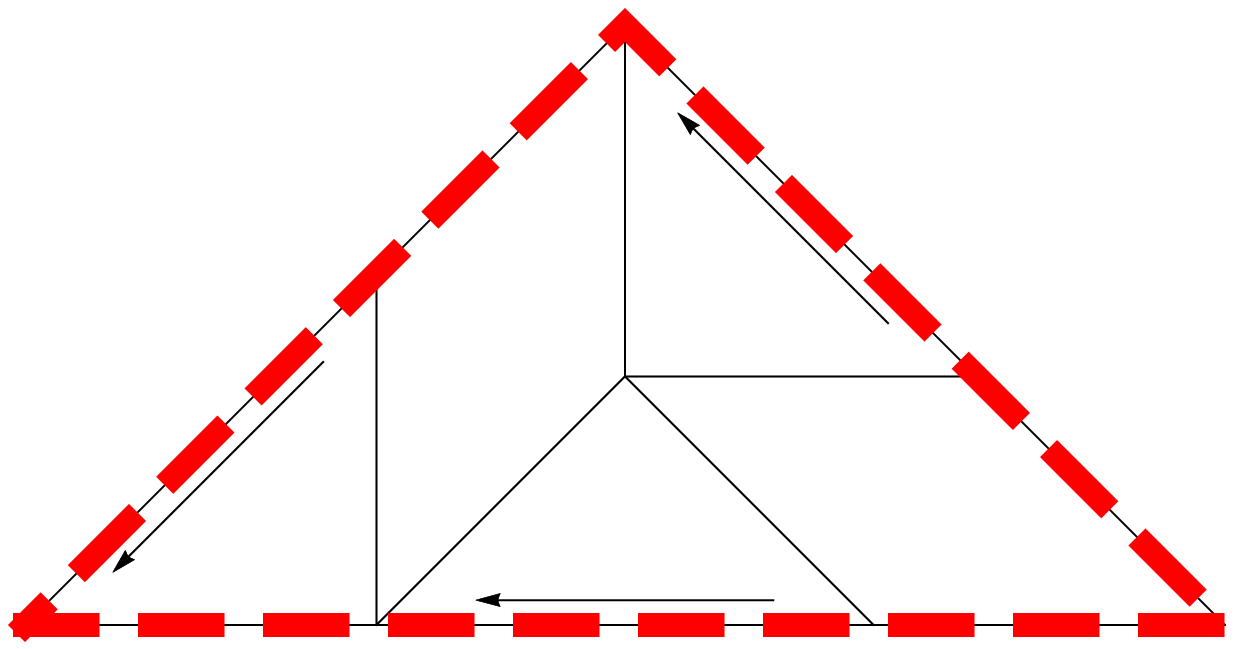}}
  {Tile A maps to 3$\times$A and 2$\times$B.}

\genfigureobject{ab_inflrule_tile_b}
  {\includegraphics[width=0.066\textwidth]{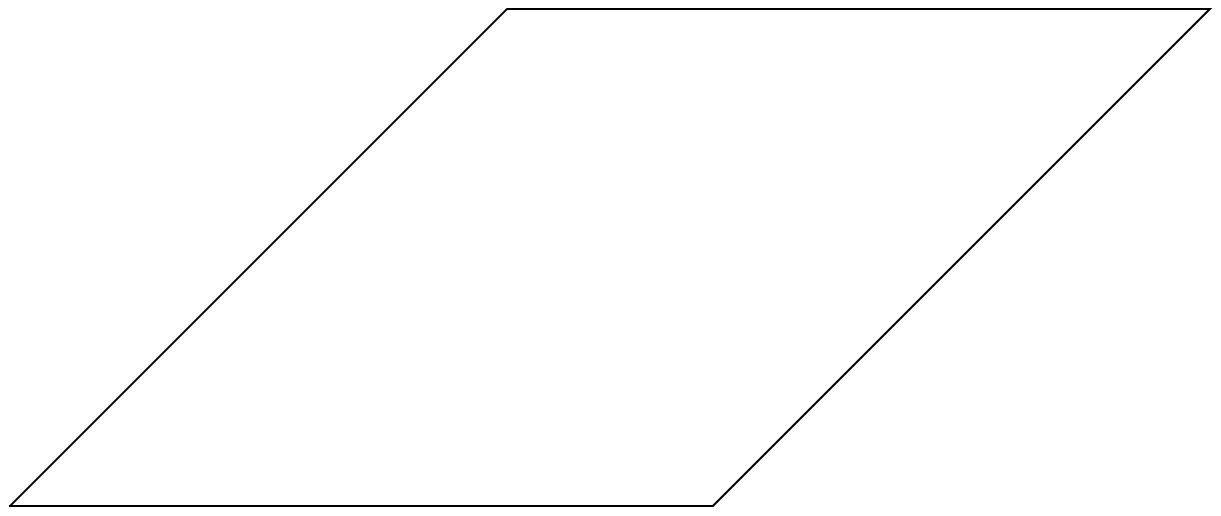} \hspace{1ex}
   \raisebox{0.01\textwidth}{$\xrightarrow{\hspace{9ex}}$} \hspace{1ex}
   \includegraphics[width=0.16\textwidth]{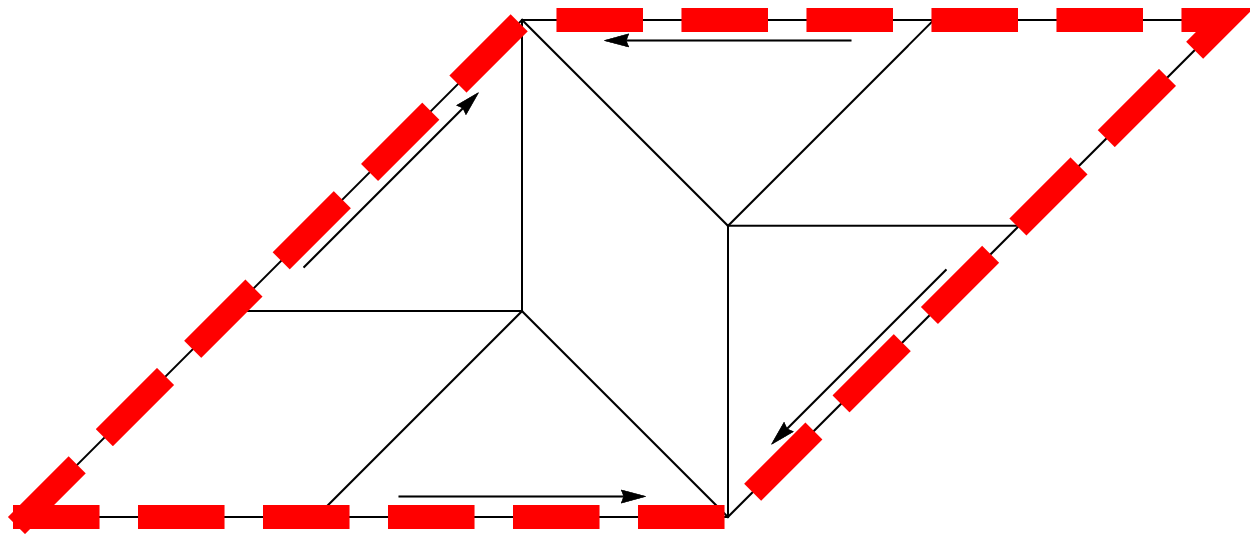}}
  {Tile B maps to 4$\times$A and 3$\times$B.}

\genfigureobject{ab_visible_tiling}
  {\includegraphics[width=0.28\textwidth]{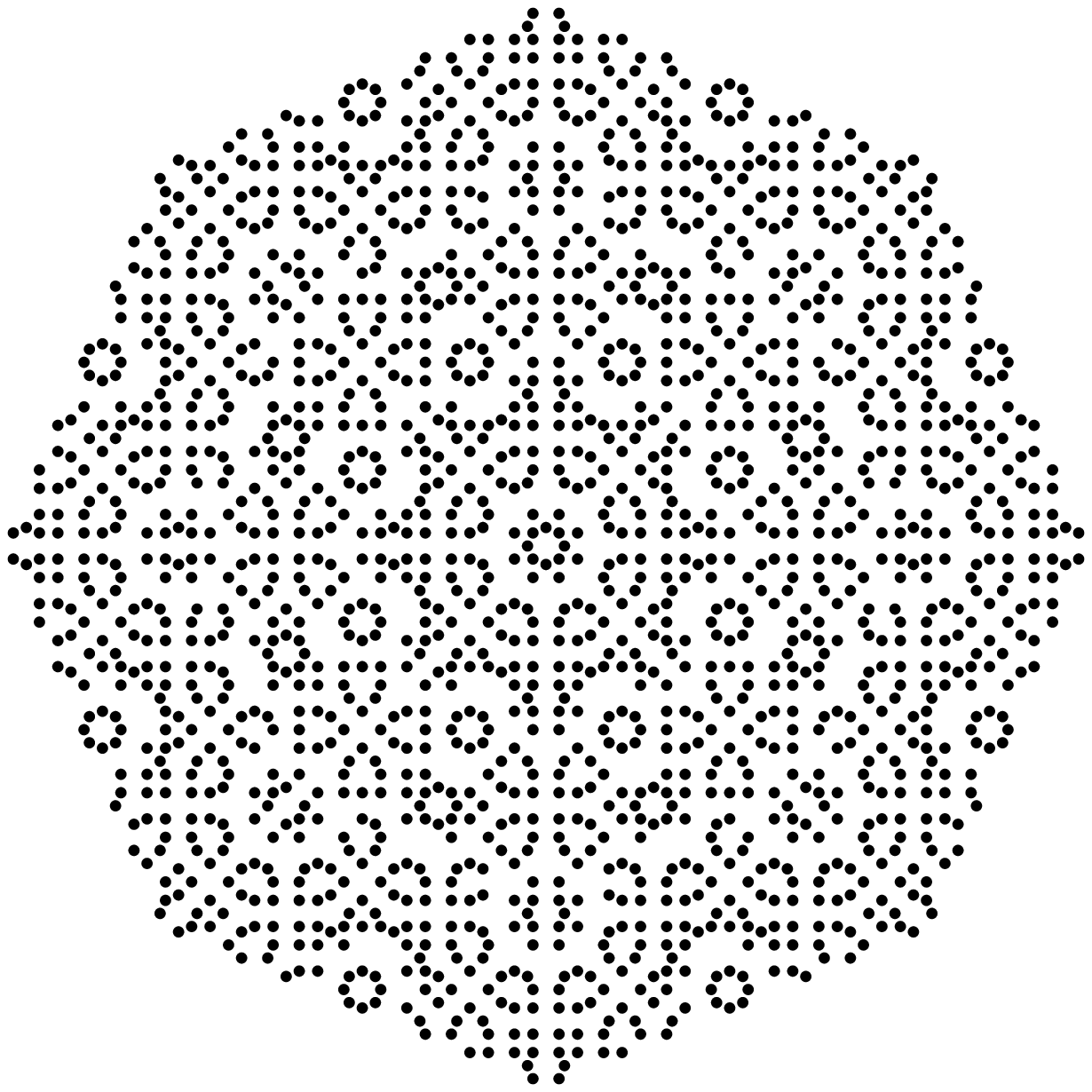} \hspace{1.8ex}
   \raisebox{0.056\textwidth}{\includegraphics[width=0.16\textwidth]{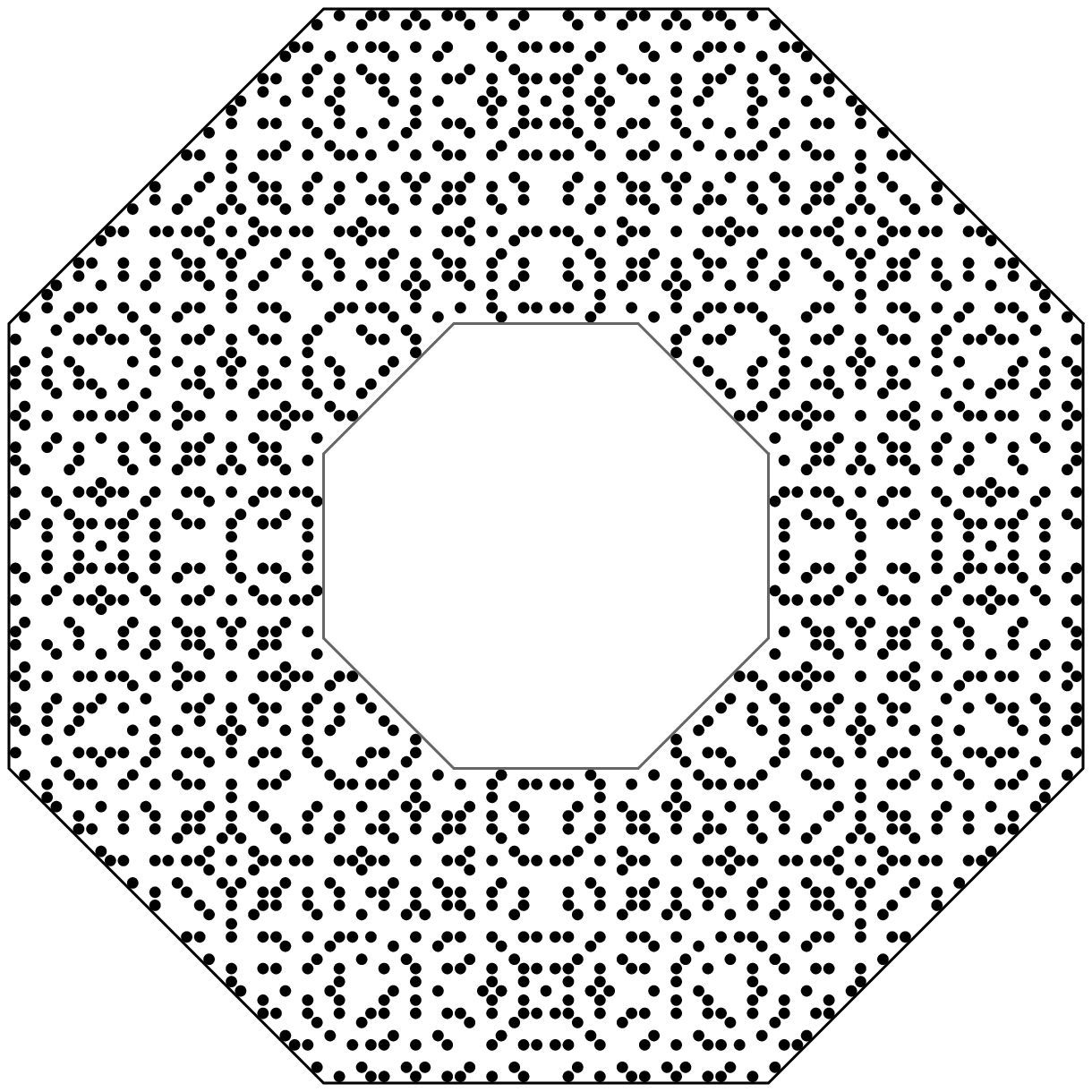}}}
  {Visible vertices of the $8$-fold symmetric \emph{Ammann--Beenker} tiling\\
   (left: direct space, right: internal space).}

\genfigureobject{ab_histogram_full}
  {\includegraphics[width=0.37\textwidth]{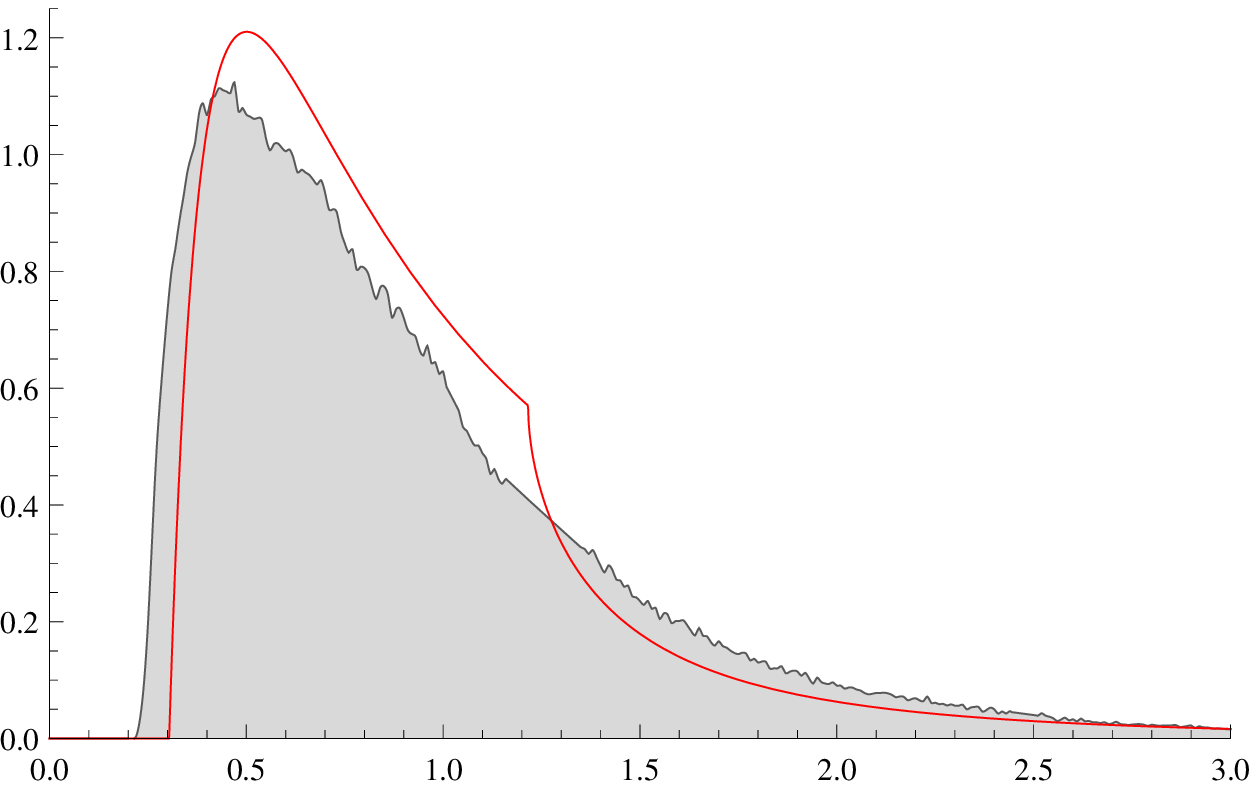}}
  {Distribution of the radial spacings of a large
   $D_8$-symmetric \textsf{AB} patch.}

\genfigureobject{ttt_circpatch_infl4}
  {\includegraphics[width=0.36\textwidth]{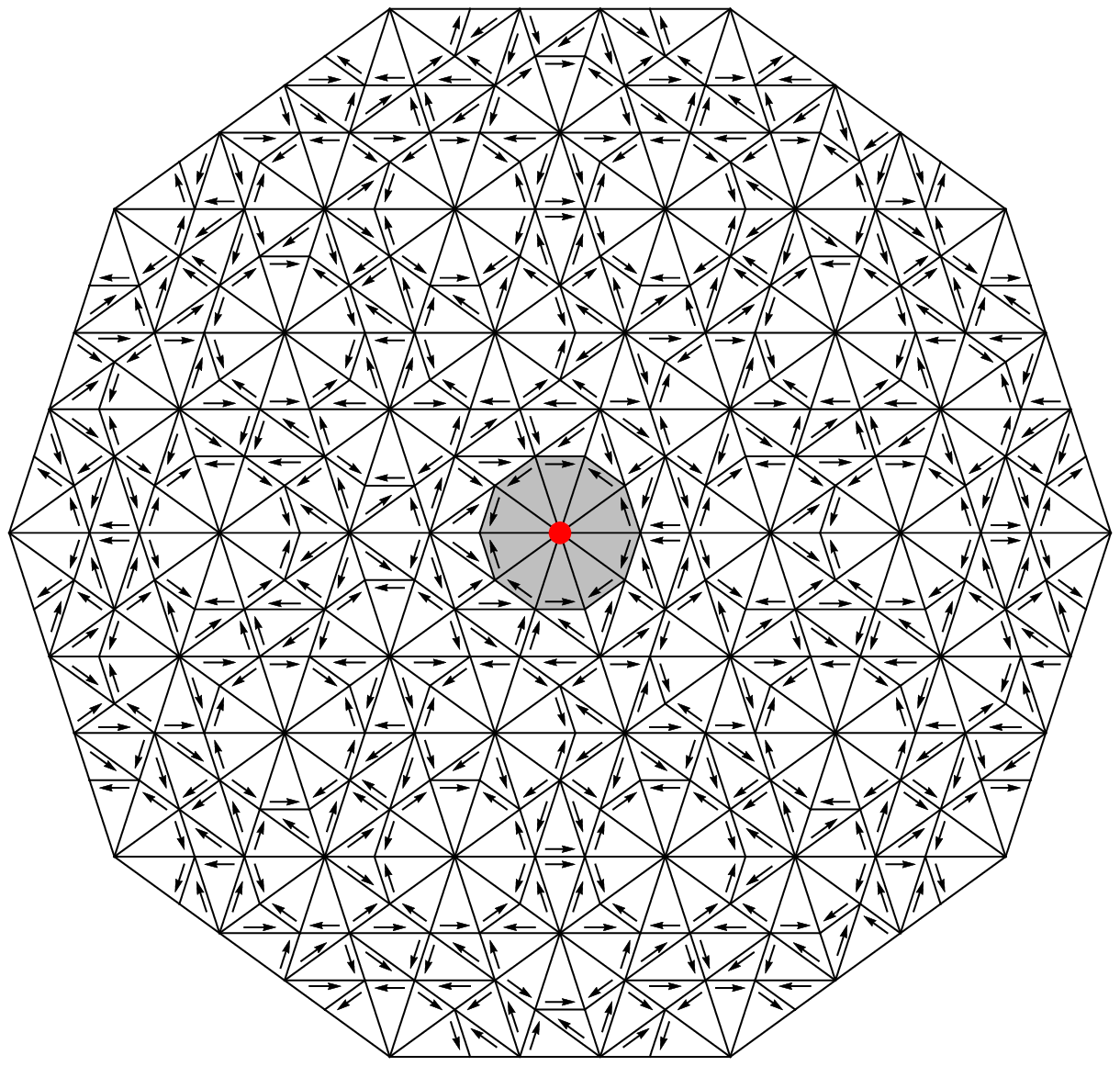}}
  {Patch of the \emph{T\"ubingen triangle} tiling
   (after 4 inflations of the central patch)}

\genfigureobject{ttt_histogram_full}
  {\includegraphics[width=0.37\textwidth]{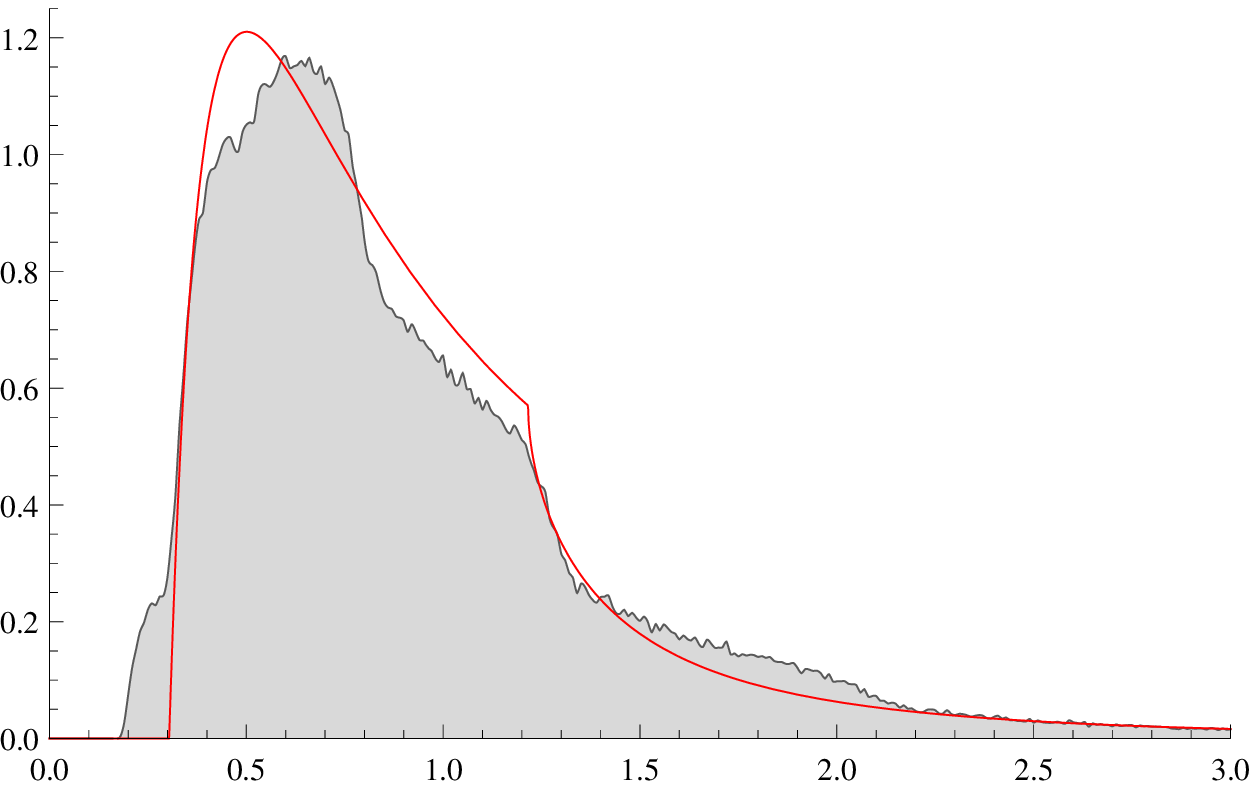}}
  {Spacing distribution of a large \emph{T\"ubingen triangle} patch.}

\genfigureobject{ttt_histogram_zoom}
  {\includegraphics[width=0.37\textwidth]{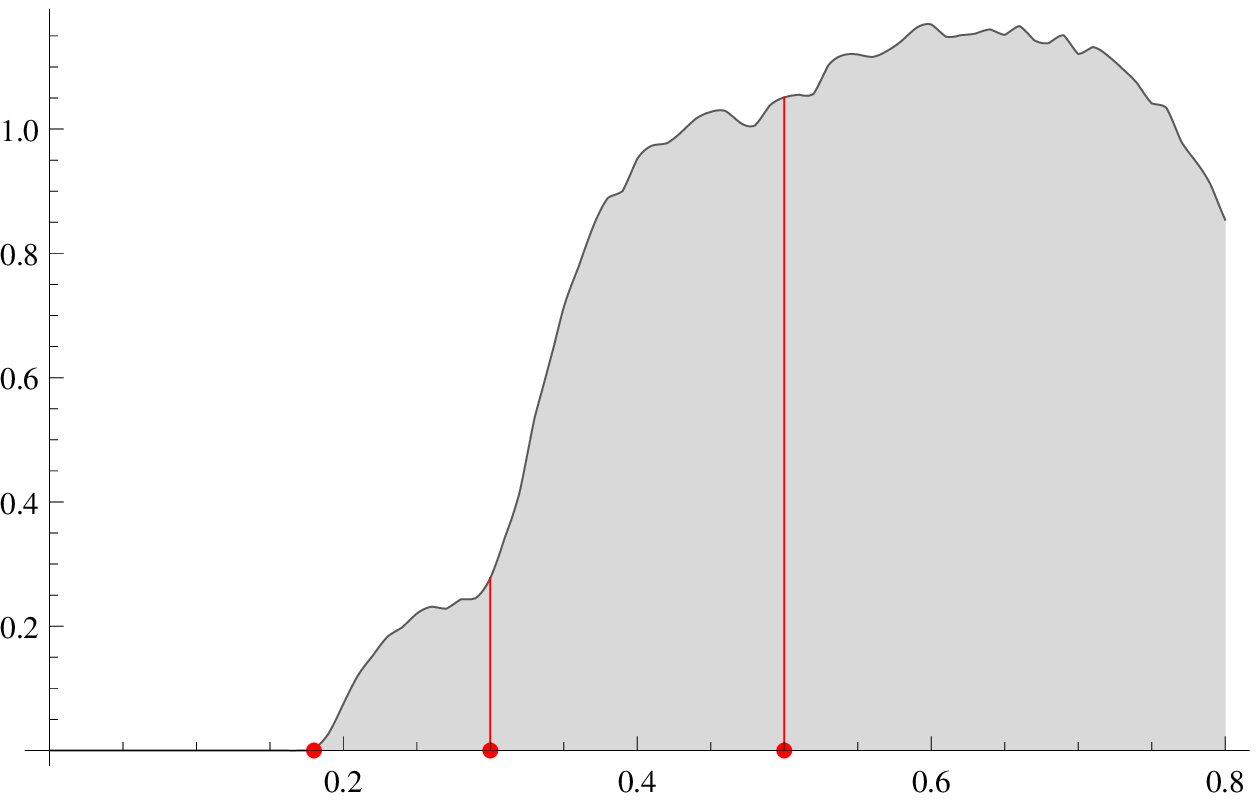}}
  {Zoom into the bulk of the \emph{T\"ubingen triangle} distribution.}

\genfigureobject{gst_visible_tiling}
  {\includegraphics[width=0.28\textwidth]{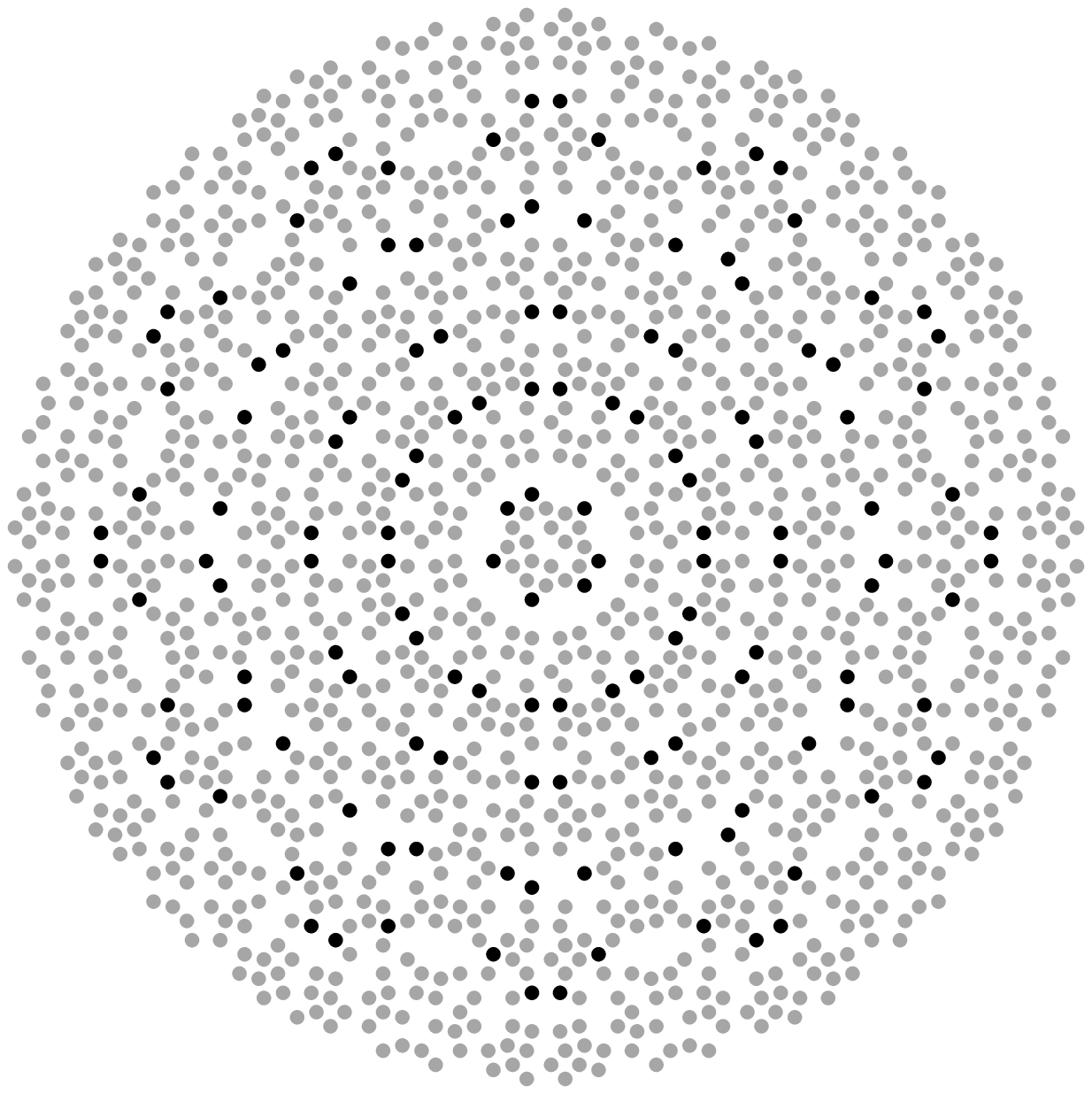} \hspace{1.8ex}
   \raisebox{0.055\textwidth}{\includegraphics[width=0.16\textwidth]{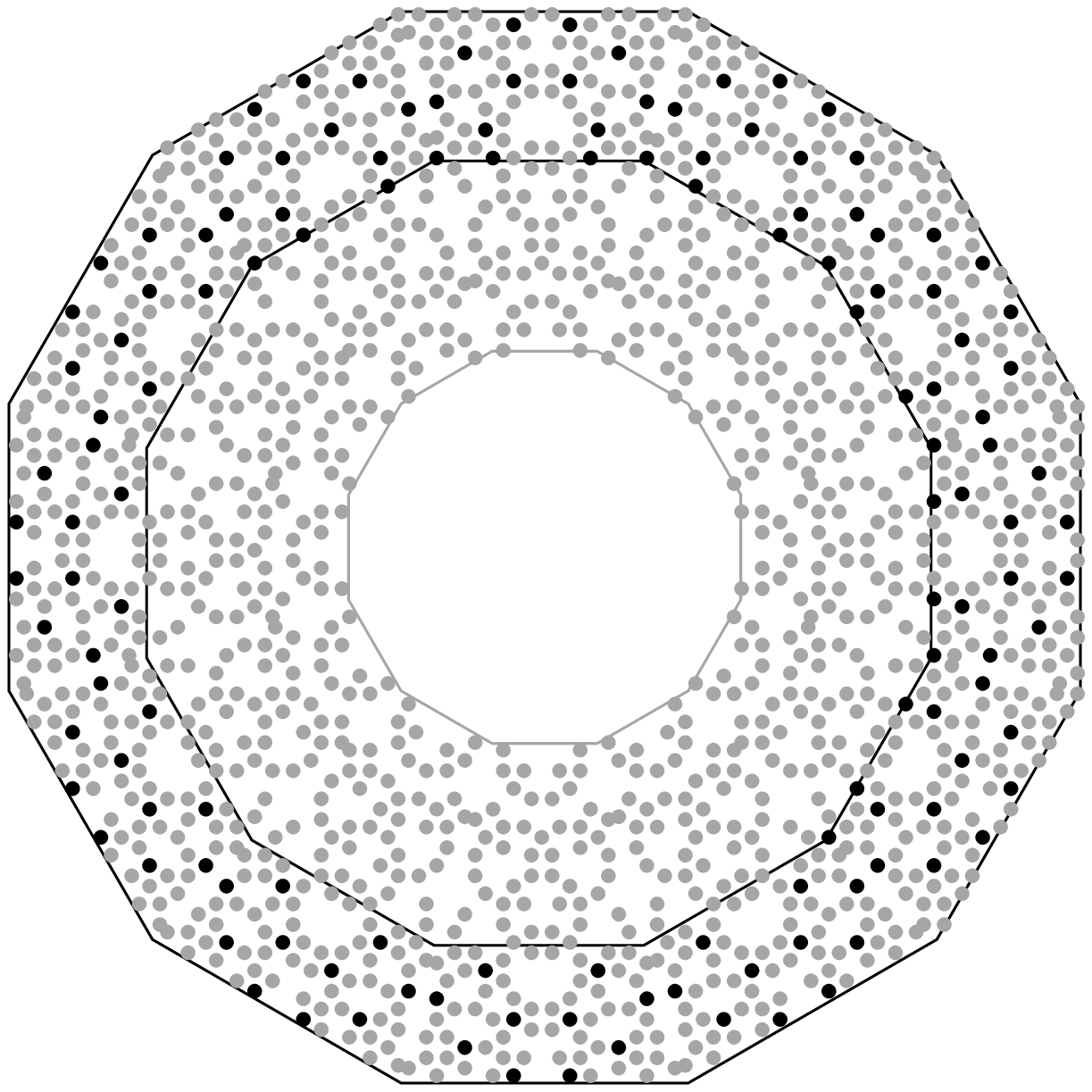}}}
  {Visible vertices of a \textsf{GS} tiling\\
   (left: direct space, right: internal space).}

\genfigureobject{gst_histogram_full}
  {\includegraphics[width=0.37\textwidth]{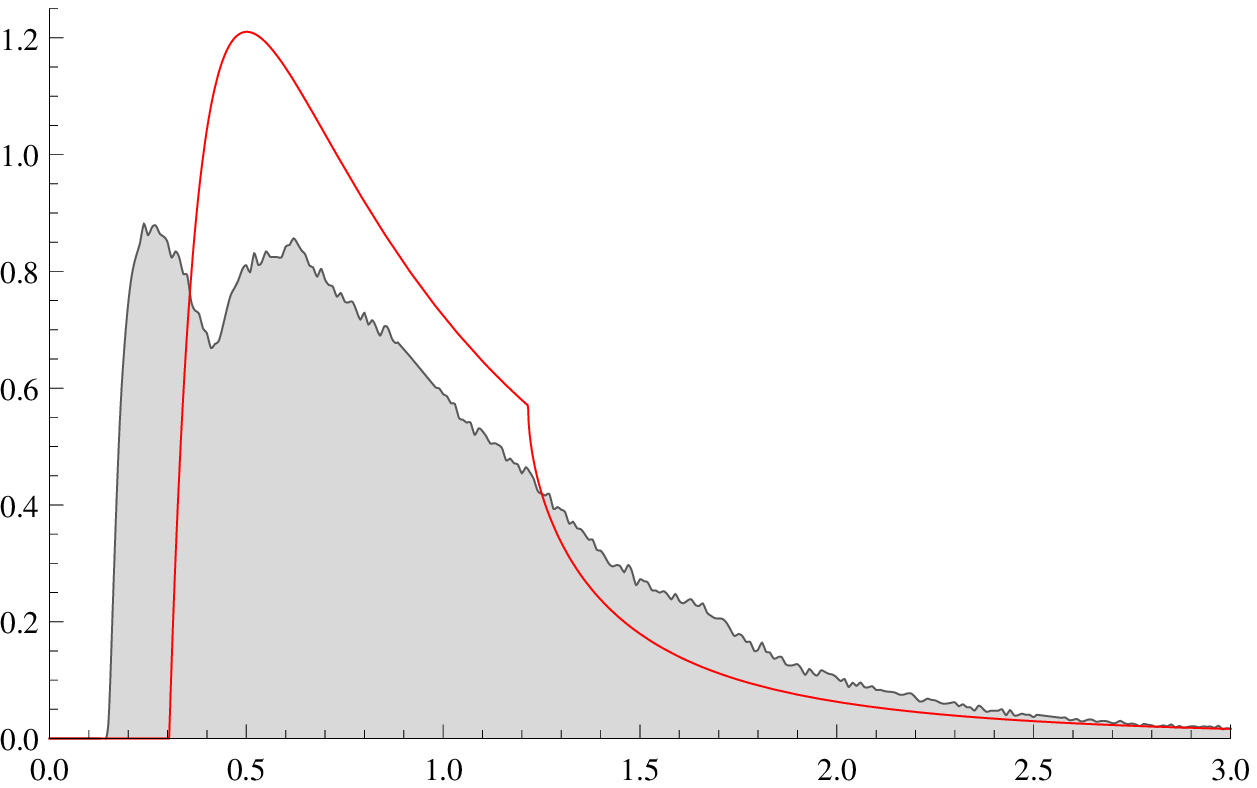}}
  {Spacing distribution of a large patch of the \emph{G\"ahler shield} tiling.}

\genfigureobject{chir_inflrule_tile_a}
  {\includegraphics[width=0.084\textwidth]{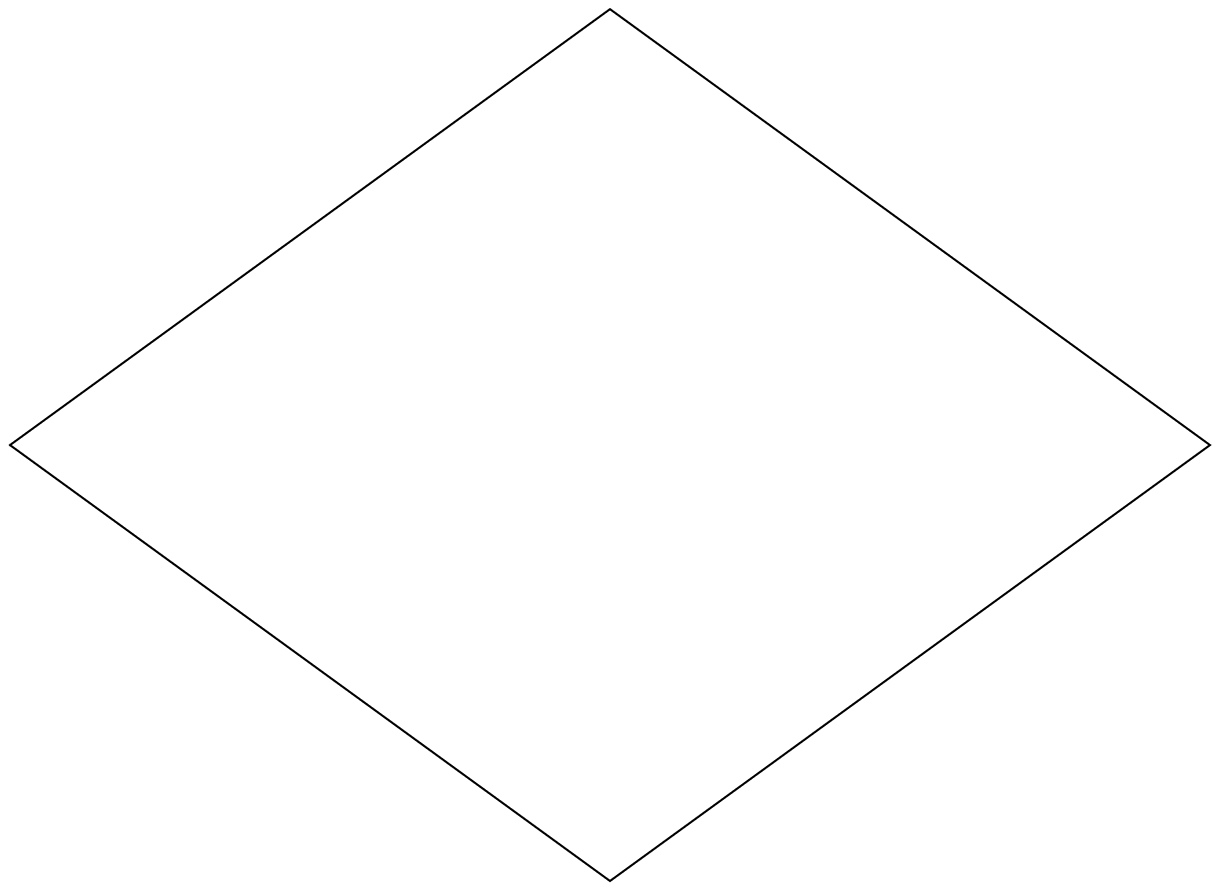} \hspace{1ex}
   \raisebox{0.026\textwidth}{$\xrightarrow{\hspace{6ex}}$} \hspace{1ex}
   \includegraphics[width=0.16\textwidth]{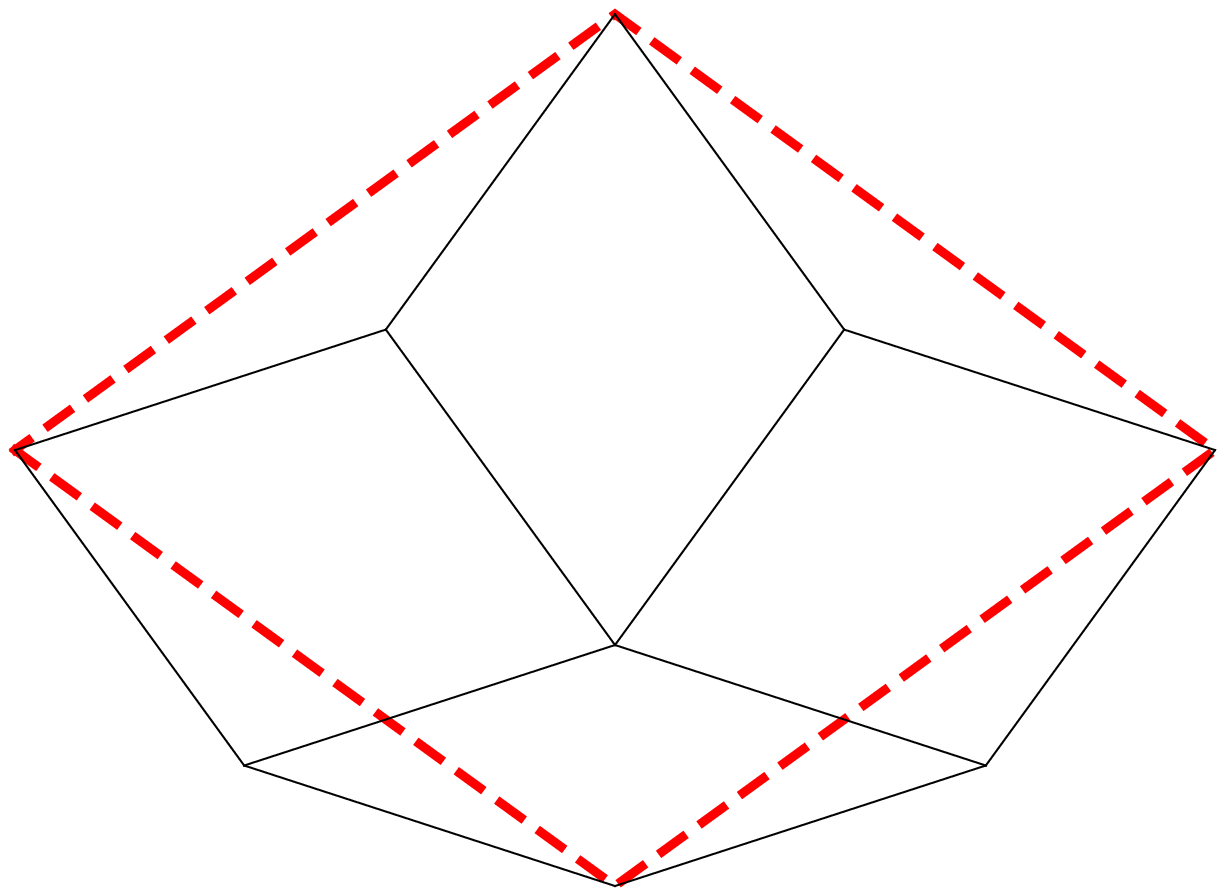}}
  {Tile A maps to 3$\times$A and 1$\times$B.}

\genfigureobject{chir_inflrule_tile_b}
  {\includegraphics[width=0.084\textwidth]{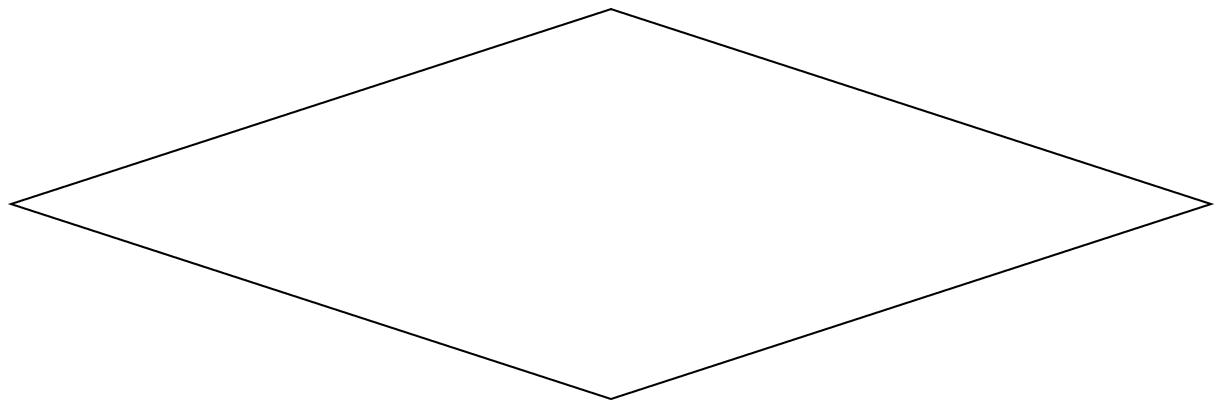} \hspace{1ex}
   \raisebox{0.011\textwidth}{$\xrightarrow{\hspace{6ex}}$} \hspace{1ex}
   \includegraphics[width=0.16\textwidth]{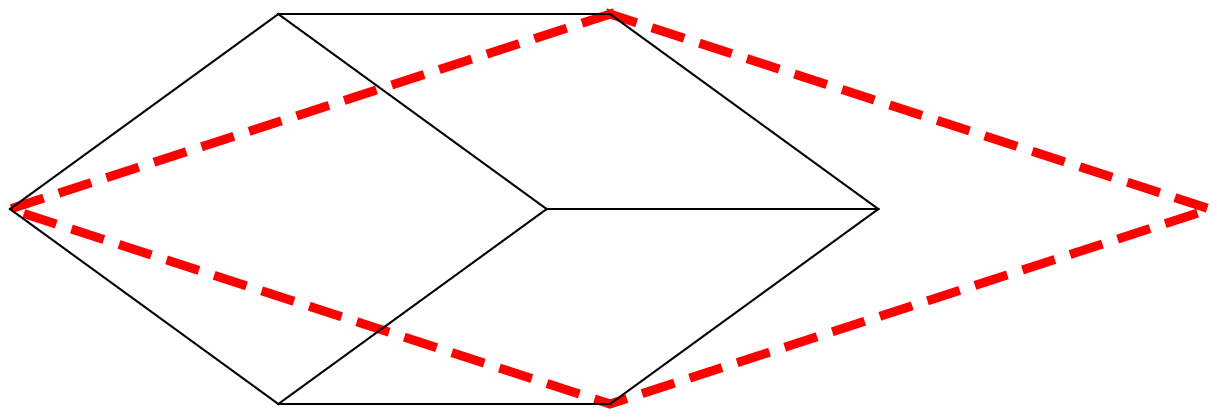}}
  {Tile B maps to 1$\times$A and 2$\times$B.}

\genfigureobject{chir_circpatch_infl4}
  {\includegraphics[width=0.36\textwidth]{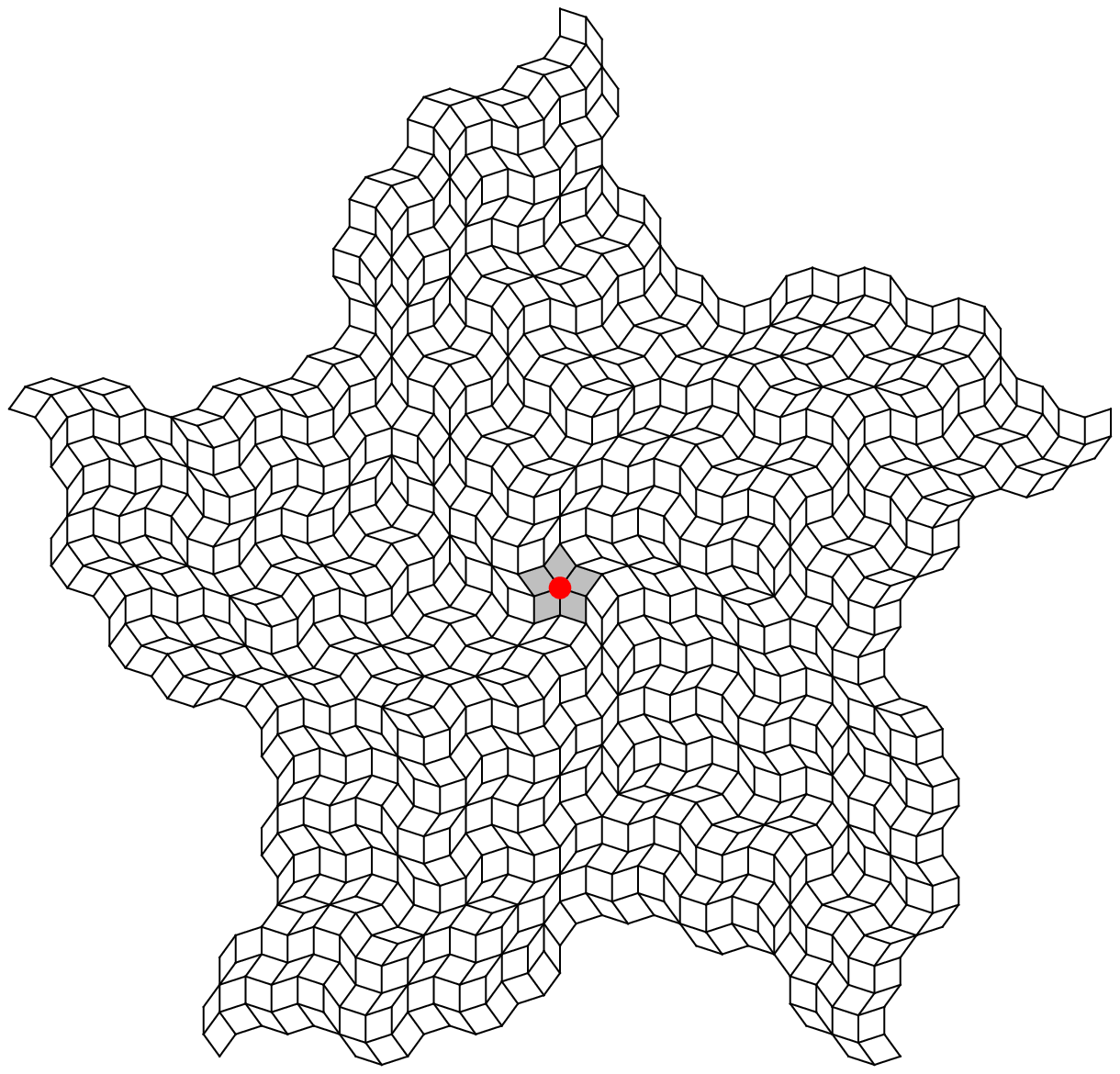}}
  {Fivefold symmetric patch of the chiral \textsf{LB} tiling
   (after 4 inflations of the \textcolor{lgray}{initial patch}).}

\genfigureobject{chir_histogram_full}
  {\includegraphics[width=0.37\textwidth]{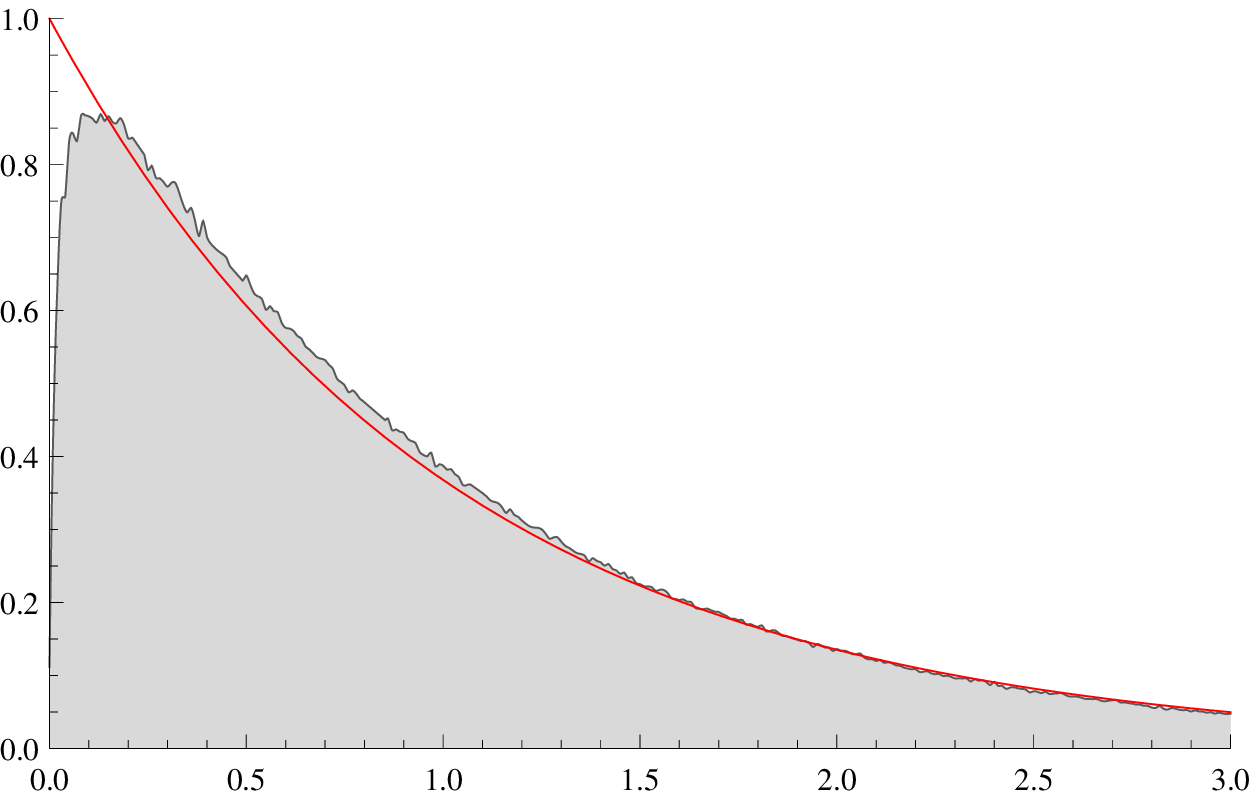}}
  {Spacing distribution of a large patch of the \emph{Lan\c{c}on--Billard} tiling.}

\genfigureobject{chair_histogram_full}
  {\includegraphics[width=0.37\textwidth]{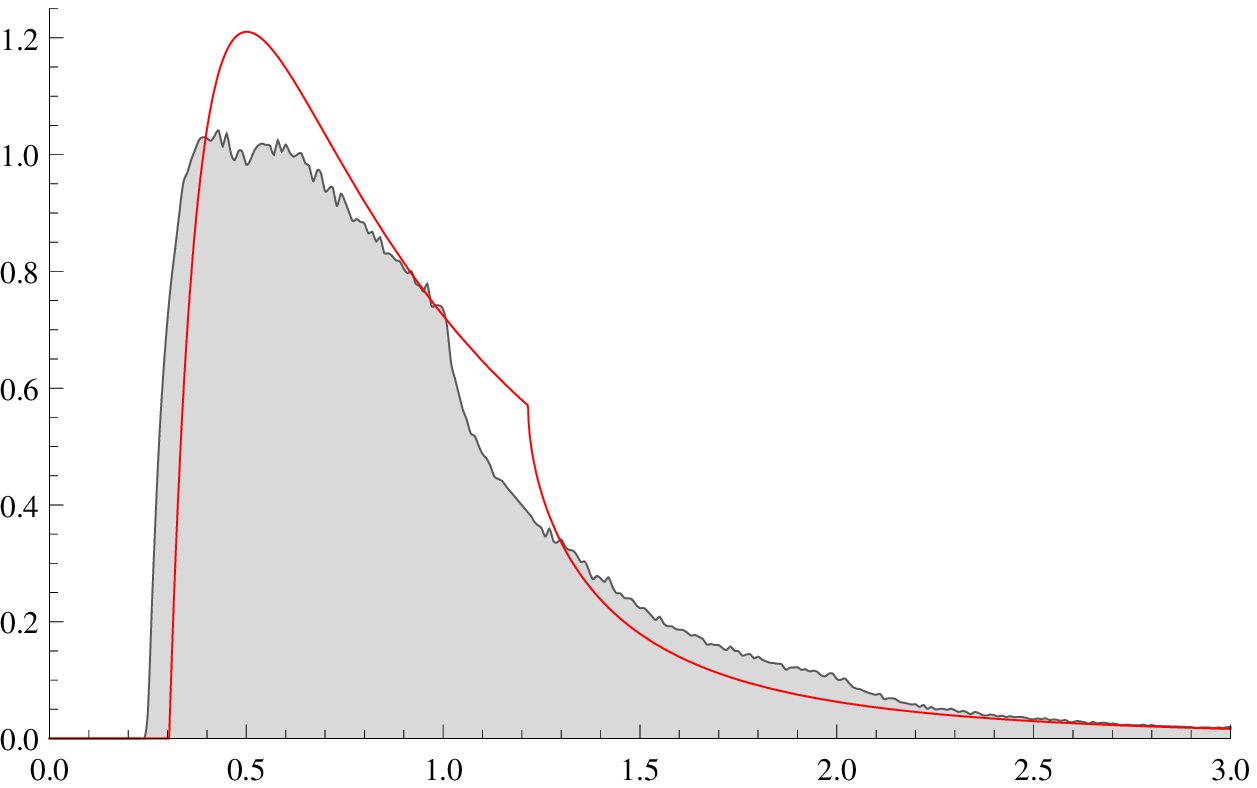}}
  {Spacing distribution of a large patch of the \emph{chair} tiling.}

\genfigureobject{badchair_histogram_full}
  {\includegraphics[width=0.37\textwidth]{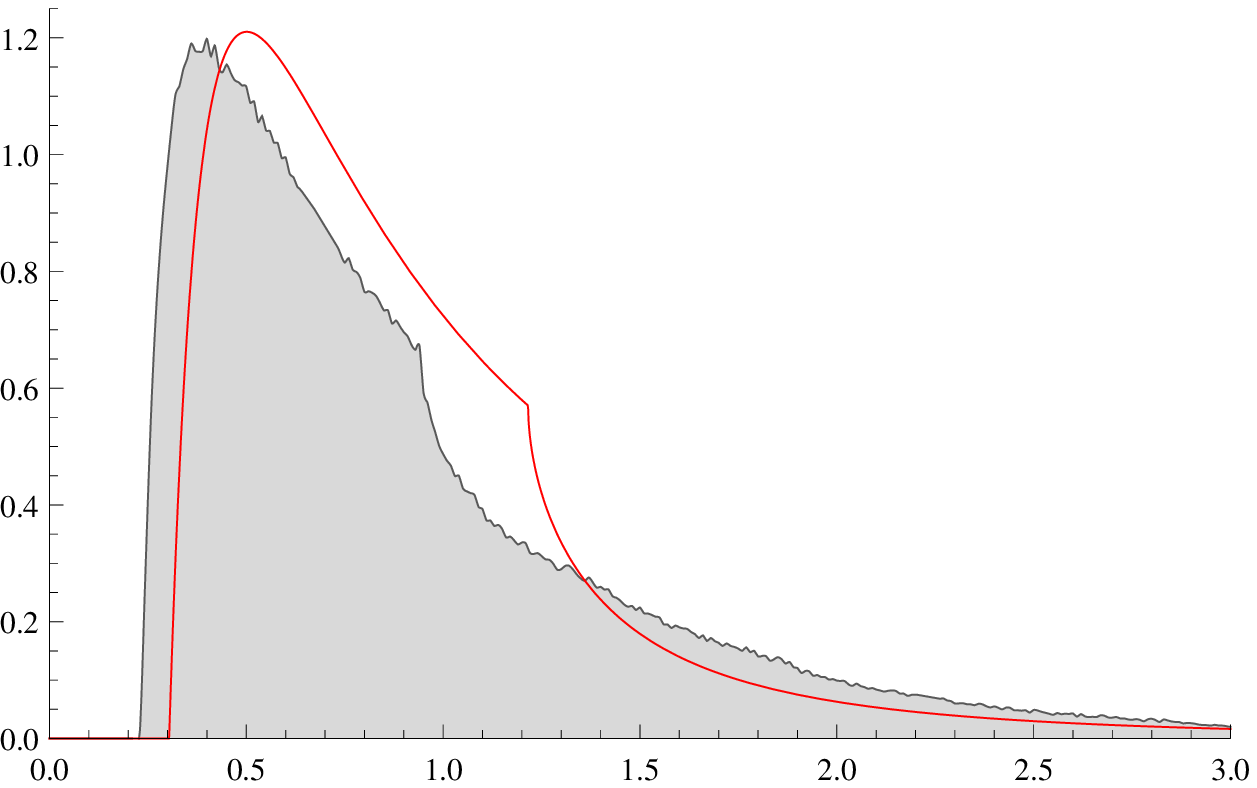}}
  {Spacing distribution of a large patch of the \emph{chair} tiling
   (using the standard gcd visibility test).}

\genfigureobject{prt_histogram_full}
  {\includegraphics[width=0.37\textwidth]{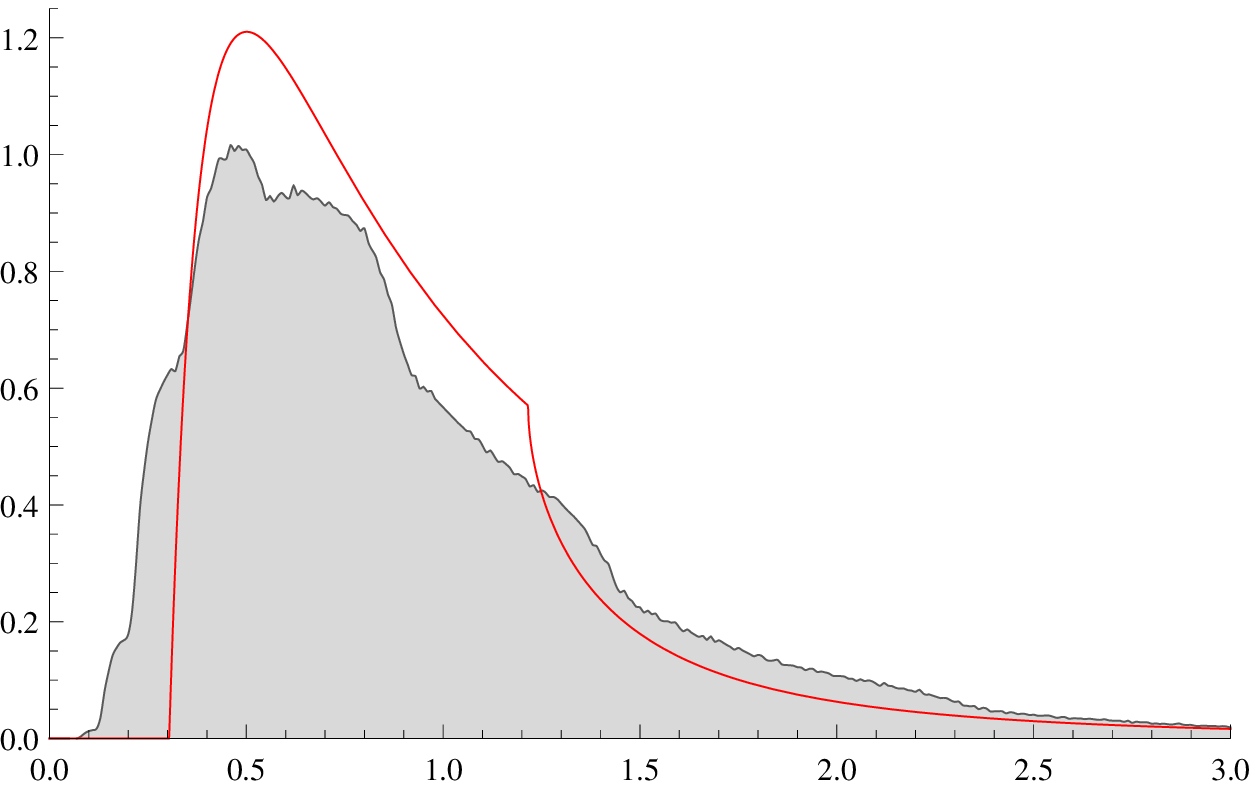}}
  {Spacing distribution of a large patch of the \emph{Penrose--Robinson} tiling.}

\genfigureobject{prt_histogram_zoom}
  {\includegraphics[width=0.37\textwidth]{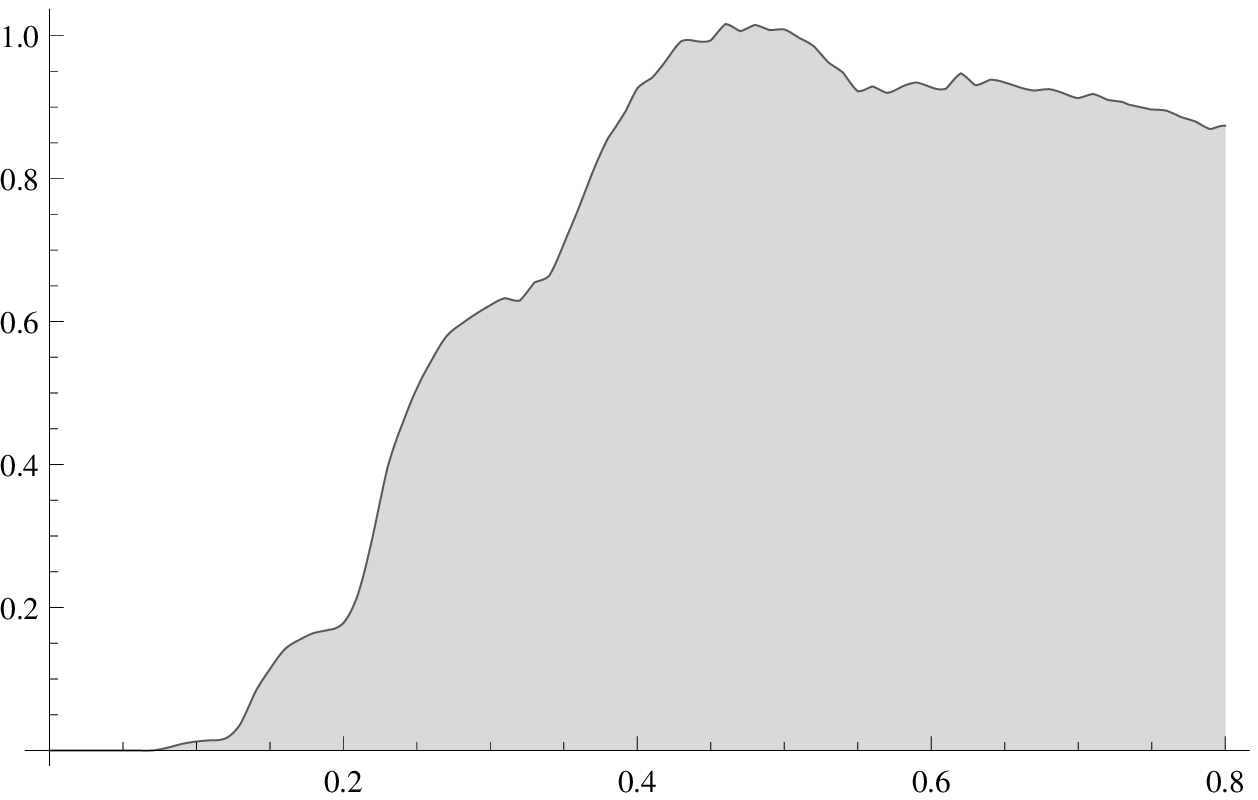}}
  {Zoom into the bulk of the \textsf{PR} distribution.}

\genfigureobject{rp_histogram_full}
  {\includegraphics[width=0.37\textwidth]{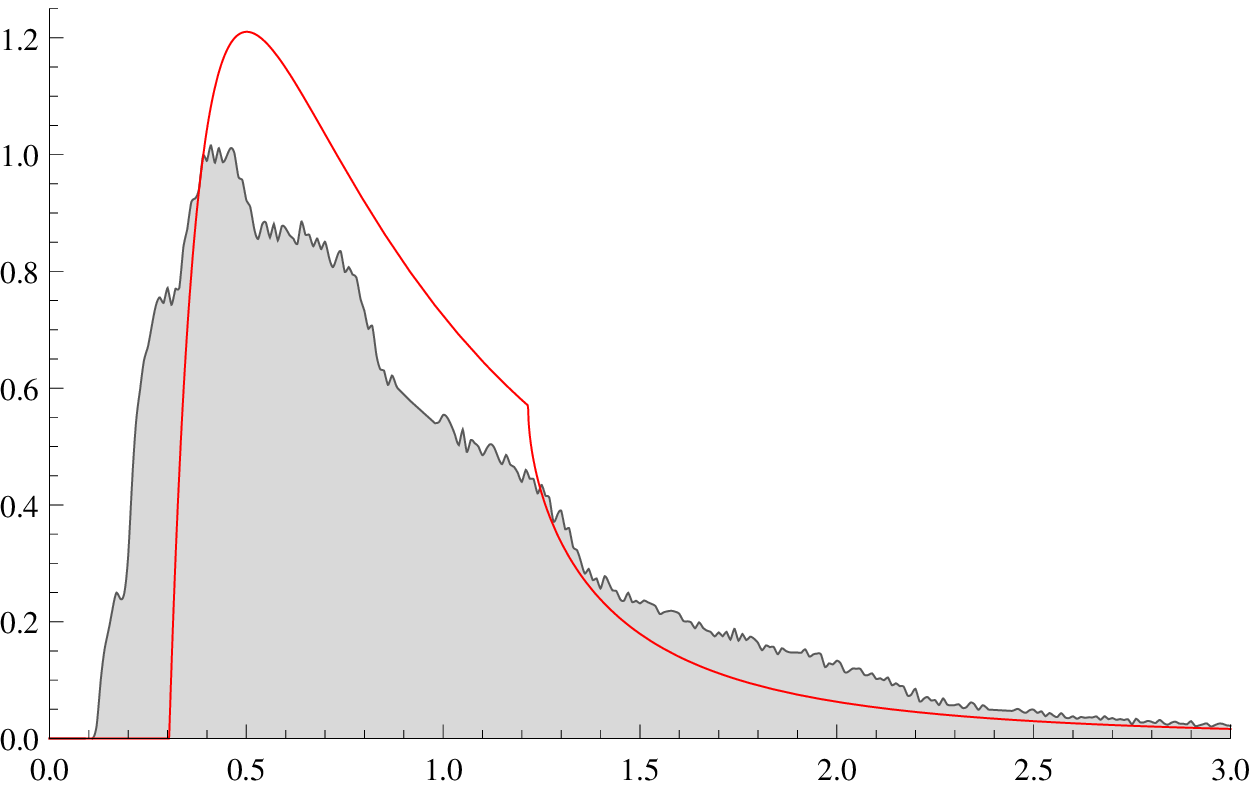}}
  {Spacing distribution of a large patch of the \emph{rhombic Penrose} tiling.}

\providecommand{\putfigure}[1]{
\def\putfigPre{insertfigure:}
\def\putfigPost{#1}
\csname \putfigPre \putfigPost \endcsname
}

\begin{document}

\title{Radial spacing distributions from planar points sets}

\author{M. Baake}
\author{F. G\"otze}
\author{C. Huck}
\author{T. Jakobi}

\address{Department of Mathematics, Bielefeld University,
         Bielefeld, Germany}
\email{\{mbaake,goetze,huck,tjakobi\}@math.uni-bielefeld.de}

\begin{abstract}
In this paper, we explore the radial projection method for locally
finite planar point sets and provide numerical examples for different
types of order. The main question is whether the method is suitable
to analyse order in a quantitive way. Our findings indicate that the
answer is affirmative. In this context, we also study local
visibility conditions for certain types of aperiodic point sets.
\end{abstract}

\maketitle

\section{Quantifying order}

When looking at physical structures, the natural question about the internal
order (of molecules, atoms, molecule clusters) arises. How to quantify order
in a good way is still largely unknown.\\
Consider a mathematical model such that the positions of the components inside the
structure are represented as a locally finite point set in $\Real^d$. We
are primarily interested in the cases ${d = 2}$ or ${d = 3}$. Let us denote the elements
of the point set as \emph{vertices}. One could now describe the order by looking
at each vertex and measure the Euclidean distance to all other vertices in the
set. This would yield a very complicated object, and comparing two such objects
resulting from different sets is going to be even more complicated.\\
This approach would be somewhat naive and also does not correspond to any physical
measurement. There are, however, methods like \emph{diffraction} (see
\cite{Hof95,Cowl95} and \cite[Ch.\,9]{BaGr13} for an introduction) that give
a lot of information about the input set. Some properties which can be analysed
by diffraction are translational repetitions and symmetries of the set.\\
Here, we present another approach, which shares some similarities with the
diffraction method, but avoids Fourier-based methods and instead works in the
direct space where the point set lives. We would like to call this the
\emph{radial projection method}, since its key ingredient is a suitable
reduction of the information coming from the point set, here implemented
by mapping a vertex to its angular component relative to some reference frame.

\section{Radial projection method}

We restrict ourselves to dimension ${d = 2}$. A possible generalisation to
higher dimensions will be discussed in Sec.~\ref{sec:conclusions}.\\
Given a locally finite point set $S \subseteq \Real^2$, we first
choose a reference point ${x_0 \in S}$. Usually, $x_0$ is chosen in such
way that it provides high symmetry (see Figure \ref{fig:ab_visible_tiling}
on page \pageref{fig:ab_visible_tiling} for an example). Now, $S$ is
thinned out by removing invisible vertices. These are the vertices that
are not observable from the reference point $x_0$, meaning that a straight line
from $x_0$ to the point, $p$ say, is already blocked by some
other point $p_0$ of the set:
\begin{equation} \exists \; p_0 \in S \; \exists \; t \in (0,1) \; :
   \; p_0 = x_0 + t \cdot (p - x_0) \; \text{.} \label{eq:occlusion_property} \end{equation}
Denote this new set of visible vertices by $V$.\\
Now, fix an ${r > 0}$ and consider the closed disk of radius $r$ around $x_0$.
Without loss of generality, we may assume $x_0 = (0,0)$. Let
$V(r)$ be the intersection of the disk and $V$. Since $S$ was chosen as locally
finite, we have $\abs{V(r)} < \infty$. We proceed by projecting each ${v \in V(r)}$ from the
reference point onto the boundary of the disk. If we write the vertex in polar
coordinates, $v = s \cdot e^{i \varphi}$ ($0 \le s \le r$), this amounts to
mapping $v$ to $\varphi$. This leaves us with a list of angles which are
then sorted in ascending order:
\[ \Phi(r) \; := \; \{ \varphi_1, \ldots, \varphi_n \} \; \text{.} \]
In fact, one has $\varphi_{i} < \varphi_{i+1}$ for all $i$ since the reduction
to visible vertices ensures that the projected vertices are distinct. The mapping
from visible vertices to their angles is therefore one-to-one.\\
By normalising with the factor $\frac{n}{2 \pi}$, the mean distance between consecutive
$\varphi_i$ becomes one. Let $d_i := \varphi_{i+1} - \varphi_i$ and define the
discrete probability measure ($\delta_x$ being the \emph{Dirac measure} at the position $x$)
\[ \nu_r \; := \; \frac{1}{n-1} \sum_{i=1}^{n-1}{\delta_{d_i}} \]
encoding these distances between consecutive angles (often denoted as discrete
\emph{spacing} distribution in the physics literature). The choice to
consider neighbouring angles is motivated by the concept of \emph{two-point correlations}
which is prominent when looking at interacting particle systems.\\
We need to know whether there exists a limit measure $\nu$,
\[ \lim_{r \rightarrow \infty}{\nu_r} = \nu \text{,} \]
in the sense of weak convergence of measures. The renormalisation step of
the angles is more a technical condition, which makes it easier to compare
$\nu_r$ for different radii $r$. It also ensures that we map the input set
to a point set of density $1$ in $\Real$. For the subsequent histograms, this
means that we measure in units of the mean distance on the $x$-axis.\\
If such a measure $\nu$ exists, we hope that it encodes enough information
about the order of the input set so that one can compare measures for
different point sets and make statements about the underlying sets.
Obviously, comparing these measures is an easier task than comparing
the original point sets.\\
Before attempting to apply this method to some interesting point sets, we
begin with some \emph{reference point sets} as limiting cases of a
potential classification.

\section{Analytic reference cases}

So far, beyond the work of \cite{ACLe13}, there are two cases which can be
fully understood analytically and which correspond to the opposite ends of
the spectrum of order. On the one end, we encounter the totally ordered
case, on the other one complete disorder.

\subsection{Perfect order / $\Integer^2$ lattice case}\label{sec:int_lattice}

Here, the choice of reference point does not matter as long as one chooses a
$x_0 \in \Integer^2$ (in \cite{MkSt10}, also a generic reference point was studied).
For simplicity, we let $x_0 := (0,0)$. A simple geometric argument then reveals
that visibility of a vertex $(x,y) \in \Integer^2$ is characterised by the
property that its Cartesian coordinates are coprime (see also
\cite{BMPl00,PlHu13}), which means $\gcd(x,y) = 1$.\\
It has long been known \cite{CoZa03} that the visible lattice points
are intimately related to the \emph{Farey fractions}
\[ \mathcal{F}_Q \; = \; \left\{ a/q: 1 \le a \le q \le Q, \;
   \gcd(a,q) = 1 \right\} \text{,} \]
here of order $Q$. Sorted in ascending order, $\mathcal{F}_Q$ is also called
a Farey series, even though it technically is a finite sequence. These
sequence are especially interesting since certain uniformity
conditions are tied to one of the most important problems in mathematics.
Denote by $\mathcal{F}_Q(i)$ the $i$th entry of the series $\mathcal{F}_Q$.
Then, the growth statement
\[ \forall \; \epsilon > 0 \; : \; \sum_{i=1}^{m}{
     \absAlt{\mathcal{F}_Q(i) - \frac{i}{m}}} = \mathcal{O}(Q^{1/2 + \epsilon}) \]
     ($m = \abs{\mathcal{F}_Q}$) is equivalent to the Riemann hypothesis
\cite{LdFr24}. Another property worth noting is the closed description of
succesive fractions, which admits enumeration formulas that make an
analytic approach possible.\\
In 2000, a proof \cite{BCZa00} was presented for the existence of a continuous limit
distribution in this case. This even holds for general star-shaped expanding regions
with some extra conditions (continuity and piecewise $\mathcal{C}^1$ for the boundary).
The density function, consisting of three regions, reads
\[ g(t) = \begin{cases}
     0, & 0 < t < \frac{3}{\pi^2}, \\
     \frac{6}{\pi^2 t^2} \cdot \log{\frac{\pi^2 t}{3}}, &
       \frac{3}{\pi^2} < t < \frac{12}{\pi^2}, \\
     \frac{12}{\pi^2 t^2} \cdot \log{\left( 2 \big{/} \! \left( 1 +
       \sqrt{1 - \frac{12}{\pi^2 t}} \right) \right)}, &
       t > \frac{12}{\pi^2},
   \end{cases}
\]
and belongs to our choice of a circular (a closed disk is placed around the reference
point $x_0$) expanding region.

\subsection{Total disorder / Poisson case}\label{sec:poisson_chaos}

On the opposite end of the spectrum, we encounter the totally disordered case. In physics
terminology, this is the realm of the ideal gas. The vertices in $\Real^2$ are
distributed according to a homogeneous spatial Poisson point process, a model also known
as \emph{complete spatial randomness} (\emph{CSR}), emphasising that points are randomly
located in the ambient space.\\
In detail, let $\mu$ denote the standard Borel-Lebesgue measure on $\Real^2$ and $V$ the
random vertex set of our ideal gas. For $A \subseteq \Real^2$, define $N(A)$ to be the
number of vertices from $V$ in $A$. Then, $V$ is characterised by the following properties:
\begin{enumerate}[label=(\hspace*{-1pt}\emph{\alph*}\hspace*{-1pt})]
  \item\label{pois_distr}{
        For each measurable $A \subseteq \Real^2$, the quantity $N(A)$
        is a Poisson random variable, which is distributed according to
        $\Pois(\lambda \mu(A))$ for a fixed ${\lambda > 0}$.}
  \item\label{pois_indep}{
        For each finite selection of disjoint $A_1, \ldots, A_k \subseteq \Real^2$,
        the quantities $N(A_1), \ldots, N(A_k)$ are independent random variables.}
\end{enumerate}
The Poisson property \ref{pois_distr} implies a condition for \emph{overlapping} vertices,
\begin{equation}
  \label{eq:pois_no_overlap}
  \lim_{\mu(A) \rightarrow 0}{\frac{\Prob(N(A) \ge 1)}{\Prob(N(A) = 1)}} = 1 \; \text{.}
\end{equation}
The probability to find more than one vertex in a volume $A$ therefore vanishes
when $\mu(A)$ goes to zero.

Fix a radius ${r > 0}$ and project the vertices from $V \cap \overline{B_r(0)}$
(the choice of reference point is arbitrary) onto the boundary $\partial{{B_r(0)}}$.
First of all, the overlapping property ensures that almost surely no overlaps occur
even after the projection.\\
Define for $\varphi_1, \varphi_2 \in [0, 2 \pi)$ with $\varphi_1 < \varphi_2$
the sector
\[ S_{\varphi_1, \varphi_2}(r) \; := \; \{ z = s \cdot e^{i \theta} \; : \;
   0 \le s \le r, \; \varphi_1 \le \theta \le \varphi_2 \} \]
between the angles $\varphi_1$ and $\varphi_2$. Let $\varphi \in [0, 2 \pi)$ be
fixed, set $\varphi_1 := \varphi, \varphi_2 := \varphi + \epsilon$ and
consider the limit ${\epsilon \rightarrow 0}$.
Since $\mu(S_{\varphi_1, \varphi_2}(r)) \rightarrow 0$, the property in Eq.~\eqref{eq:pois_no_overlap}
implies that there is at most one projected vertex at the location $\varphi$.\\
Now, select a subinterval $[a, b]$ of $[0, 2 \pi]$ and study the amount $N(a, b)$
of projected points inside $[a, b]$. The vertex count in the sector $S_{a, b}(r)$ completely
determines the quantity $N(a, b)$, which, by using property \ref{pois_distr}, is a
Poisson random variable with intensity
$\lambda \mu(S_{a, b}) = \lambda \frac{r^2 \ell}{2}$ (with $\ell := \abs{a - b}$
the length of the interval),
\[ N(a, b) \; \sim \; \Pois(\lambda \frac{r^2 \ell}{2}) \; \text{.} \]
The mean number of points in $B_r(0)$ is $\lambda \pi r^2$. Normalising the angles with
$\frac{n}{2 \pi}$, $n$ the number of vertices inside $B_r(0)$, generates a new CSR with
intensity ${\lambda = 1}$ on $\Real_{+}$ in the limit ${r \rightarrow \infty}$. The
independence property \ref{pois_indep} carries over to dimension one in an analogous way.\\
The distance between consecutive points of a spatial Poisson process in $\Real$ is known
to be exponentially distributed with density function
\[ f_{\lambda}(t) \; = \; \begin{cases}
     \lambda \exp(-\lambda t), & t \ge 0,\\
     0,                        & t < 0 \; \text{.}
   \end{cases}
\]
In the probabilistic (temporal) interpretation of a Poisson process, this is the
distribution of the waiting time between jumps. Our reference densities
therefore have these shapes:
\putfigure{z2lattice_vs_poisson}
The graphs in Figure \ref{fig:z2lattice_vs_poisson} were produced by numerical
evaluation, using $N \approx 1.98 \cdot 10^6$ angles in the $\Integer^2$ lattice
case (radius ${r = 2900}$), and $N \approx 1.96 \cdot 10^6$ angles in the Poisson
case. The analytic density functions perfectly match the graphs, which gives a
hint at how large the amount of samples has to be in general to produce appropriate
approximations.\\
Our interest now is to study other point sets and to check how they fit into this
picture. Can one expect some kind of interpolation behaviour between the two
reference densities? The primary focus will be on vertex sets coming from aperiodic
tilings, since these feature both a repetitive structure but also disorder. In terms
of density functions, one might then expect some ``mixture'' of the $\Integer^2$ and
the Poisson case.\\
We point out that the existence of a limit distribution is known in the two reference
cases. In all other considered cases, we assume that the distribution exists, which
is plausible from the numerics. A first step to prove this is given in
\cite[Thm.\,A.1]{MkSt13}. Since the release of the article's preprint, further results
\cite{MkSt14} became available, wherefore we now know the existence of the distribution
for regular model sets.

\section{Numerical approach}\label{sec:num_approach}
As mentioned in Sec.~\ref{sec:int_lattice}, the analytic approach for the
integer lattice case is based on the theory of Farey fractions. This
framework does not extend properly to arbitrary locally finite point sets.
And even for subsets of $\Integer$-modules (like all our covered examples
are), this fails since the key property, the closed description for neighbouring
fractions mentioned in \ref{sec:int_lattice}, does not hold anymore -- or
at least not in an obvious way. One would first need to extend the
notion of Farey fraction in a well-defined manner to $\Integer$-modules, but
even then it is still unclear whether the approach presented in \cite{BCZa00}
carries over.\\
From this perspective, an initial approach through numerical methods was chosen. The basic
idea is to generate a large list of vertices such that the list needs only a minimal amount
of trimming to have a circular shape. Since our focus is on aperiodic tilings, the primary
step consisted in creating large patches of these, from which we could then extract the
vertex sets with the required properties. The trimming is unavoidable since both feasible
methods introduce restrictions on the shape of the generated patch.\\
There are essentially three methods to produce aperiodic tilings of the plane. The first
one is by defining a set of prototiles with matching rules. This method is not suitable for
the purpose of implementation. We therefore focus on the alternatives, namely inflation and
projection.

\subsection{Inflation rules}\label{sec:inflation_method}

Probably the most prominent method is via inflation of prototiles. For example, the
\emph{T\"ubingen triangle} tiling (abbreviated as \textsf{TT}) is produced from two
prototiles \cite[Ch.\,6.2]{BaGr13}, both with edge length ratio ${\tau : 1}$. Here,
$\tau$ is the golden mean, which also serves as the inflation factor. The first
tile, denoted as \emph{type A}, is inflated according to the scheme shown in
Fig.~\ref{fig:ttt_inflrule_tile_a} (rescaled version indicated in \textcolor{red}{red})
\putfigure{ttt_inflrule_tile_a}
while \emph{type B} follows the rule shown in Fig.~\ref{fig:ttt_inflrule_tile_b}.
\putfigure{ttt_inflrule_tile_b}
One can see from the rules that the prototiles appear in both chiralities in the resulting
tiling. The reflected tiles are simply inflated via the reflected rules.\\
It can be shown that, for properly chosen edge lengths, the resulting vertex set
lives in $\CyI{5}$ with ${\zeta_n := \exp(2 \pi i / n)}$ a primitive $n$-th root
of unity. The first step, however, is to generate the tiling patch itself and
afterwards to extract the vertices. We start with one of the prototiles and apply
the inflation rule a few times, inspecting the result for symmetric subpatches in
each step. In this case, the inflation rule applied to one prototile of type
\emph{A} produces the patch shown in Fig.~\ref{fig:ttt_apatch_infl5} (subpatch
shaded in \textcolor{lgray}{grey}):
\putfigure{ttt_apatch_infl5}
Now, one can isolate the indicated subpatch and use it as initial patch for the
inflation. From the computational point of view, this imposes some difficulties.
We formulate these for general modules $\CyI{n}$, while keeping in mind the example of the
\textsf{TT} tiling (${n = 5}$) for illustrative purposes.
\begin{enumerate}[label=(\hspace*{-1pt}\arabic*\hspace*{-1pt})]
  \item\label{compute_integerarith}{
        Inflation steps are applied iteratively. This quickly leads to accumulation
        of numerical errors. To avoid this, we solely employ integer arithmetic and
        only switch to floating-point when computing the angular component $\arctan(y/x)$
        of a vertex $(x,y)$.}
  \item\label{compute_elemencoding}{
        Elements of $\CyI{n}$ need to be encoded exactly. These types of $\Integer$-modules
        can be written as
        \[ \hspace{8ex} \CyI{n} \; = \; \{ a_0 + a_1 \zeta_n + \ldots +
           a_{r-1} \zeta_n^{r-1} \; : \; a_i \in \Integer \} \]
        ($r = \phi(n)$ the \emph{Euler totient} function) and therefore only require $r$
        integers to encode one element (resulting in a vertex size of $4 \times 4 = 16$
        bytes for the \textsf{TT} if one uses standard $32$-bit integers). The vertex byte
        count is in fact significant, see point \ref{compute_visibility}.}
  \item\label{compute_tileencoding}{
        The inflation rule applies to prototile objects, so we have to keep a tile list
        during the patch construction. Because of \ref{compute_integerarith}, we want an
        exact encoding for list elements. We represent a tile using the type (\emph{A}/\emph{B}
        for \textsf{TT}), the chirality (not always needed), a reference point of the tile
        (exact in the $\Integer$-module case, see \ref{compute_elemencoding} above) and a
        rotation of the tile around the reference point. This requires a quantisable angle
        (the tile is only allowed to appear with a finite number of distinct rotations),
        which fails when one considers for example the famous \emph{pinwheel} tiling \cite{Radi99}.
        \begin{table}[H]
          \caption{Prototile bit encoding for the \textsf{TT} tiling.}
          \begin{tabular}{c|c|c}\toprule
            property  & states            & bit count \\
            \midrule
            type      & A / B              & $1$\\
            chirality & normal / mirrored  & $1$\\
            reference & --                 & $4 \cdot 32$ ($ = 16$ bytes)\\
            rotation  & $\{0, \ldots, 9\}$ & $4$\\
            \bottomrule
          \end{tabular}
          \label{tab:prototile_bit_encoding}
        \end{table}}
  \item\label{compute_decomposition}{
        The prototile description is only helpful while growing
        the patch, but becomes cumbersome as soon as one is
        interested in raw vertex data. Each prototile
        object decomposes into a bunch of vertices (three for
        the \textsf{TT}). Applying a decomposition step to each
        prototile in the output list yields a list with many
        duplicate vertices, requiring an additional step to reduce the
        list to unique vertices. This involves constantly accessing
        the list to locate already present vertices, making it
        preferrable to have a low element byte count.}
  \item\label{compute_visibility}{
        The determination of visibility of a single vertex is
        generally very different from the $\Integer^2$-case, where
        the test consisted of computing the $\gcd$ of the two
        coordinates. In the generic case, we have to consider the
        whole set of unique vertices to determine the visibility of
        one vertex by doing a geometric ray test (see Eq.
        \eqref{eq:occlusion_property}). It proved to be
        more efficient to combine the removal pass for
        unique vertices with the visibility test pass and to use
        custom data structures to further speed up the process.}
\end{enumerate}

The computation time mentioned in \ref{compute_visibility}, which is $\mathcal{O}(n)$,
is not to be underestimated ($n$ being the total amount of vertices collected
at some point), and led to the investigation of cases with tests having similar
complexity as $\Integer^2$, which is just $\mathcal{O}(1)$.\\
To summarise, there are roughly three steps: Growing a large circular
patch, removal of duplicate vertices together with the visibility test, and
finally mapping vertices to angles followed by proper normalisation.\\
A simple optimisation consists of removing redundancy imposed by
symmetry of the input set. For example, the $\gcd$ is fixed under sign changes
of the parameters. It also is $D_4$-symmetric, wherefore it suffices to
consider the halved upper-right quadrant of the $\Integer^2$ lattice.

\subsection{Model set description / cut-and-project}\label{sec:modelset_method}

A different method for constructing tilings is given by the cut-and-project method. The
advantage is that it directly yields vertices of the tiling and does not require keeping
track of the adjacency information. Another reason for choosing this description, if
applicable, is that some configurations admit a much easier condition to determine
visibility of a given vertex by using local information only. In this regard, such
cases are very similar to $\Integer^2$ together with the $\gcd$-test.\\
In a simplified setting, let $(\Real^d, \Real^k, \mathcal{L})$ be a
triple and $\pi, \intproj$ projections satisfying the following
conditions:
\begin{enumerate}[label=(\roman*),align=right]
  \item{$\mathcal{L}$ is a lattice in $\rprod{d}{k}$;}
  \item{$\pi : \rprod{d}{k} \rightarrow \Real^d$,
        with $\pi|_{\mathcal{L}}$ injective;}
  \item{$\intproj : \rprod{d}{k} \rightarrow \Real^k$,
        with $\intproj(\mathcal{L}) \subset \Real^k$ dense.}
\end{enumerate}
This setup is called a \emph{cut-and-project scheme} (\emph{CPS}). If we define
$L := \pi(\mathcal{L})$, the conditions above induce 
${\star : L \rightarrow \Real^k}$, the \emph{star map}. The lattice can then be written
as $\mathcal{L} = {\{ (x, x^{\star}) : x \in L \}}$ and one usually encodes the
CPS in a diagram. The right hand side in Figure \ref{fig:general_cut_and_project}
describes the \emph{internal space}, the left one the \emph{physical space} (since this
is where the point set of the tiling itself lives).\\
\begin{figure}[H]
  \centering
  \begin{tikzpicture}[description/.style={fill=white,inner sep=2pt}]
    \matrix (m) [matrix of math nodes, row sep=2em,
                 column sep=2.5em, text height=1.5ex,
                 text depth=0.25ex]
    { \Real^d            \mnc \rprod{d}{k}     \mnc \Real^k \\
      \pi(\mathcal{L})   \mnc \mathcal{L}      \mnc \intproj(\mathcal{L}) \\
      L                  \mnc                  \mnc L^{\star}  \\ };
    \path[->,font=\scriptsize]
    (m-1-2) edge node[auto,swap] {$\pi$} (m-1-1)
    (m-1-2) edge node[auto] {$\pi_{\operatorname{int}}$} (m-1-3)
    (m-2-2) edge node[auto,swap] {$1$-$1$} (m-2-1)
    (m-2-2) edge node[auto] {} (m-2-3)
    (m-3-1) edge node[auto] {$\star$} (m-3-3);
    \path[solid,font=\scriptsize]
    (m-2-1) edge[double,double distance=3pt] node[auto] {} (m-3-1)
    (m-2-3) edge[double,double distance=3pt] node[auto] {} (m-3-3);
    \path[->,font=\scriptsize,left hook-latex]
    (m-2-1) edge node[auto,swap] {} (m-1-1)
    (m-2-3) edge node[auto,swap] {dense} (m-1-3)
    (m-2-2) edge node[auto,swap] {} (m-1-2);
  \end{tikzpicture}
  \caption{General case of a $\Real$-CPS.}
  \label{fig:general_cut_and_project}
\end{figure}

Details about the generic definition can be found in \cite{Schl98,BaGr13}.
Given a CPS as defined above, a \emph{model set} then arises from choosing
a subset $W \subseteq \Real^k$ (with certain conditions) and considering the set 
\[ \mbox{\Large $\curlywedge$}(W) \; := \; \{ x \in L \; : \; x^{\star} \in W \} \text{.} \]
The subset $W$ is called the \emph{window} of the model set (also denoted
as acceptance region or occupation domain).
It can be shown that point sets of certain aperiodic tilings can be
generated using this description. This is also the important aspect for
our implementation purpose, since the main work now consists of generating
a suitable ``cutout'' ${L_0 \subset L}$ and then applying the window
condition ${x^{\star} \in W}$ to each ${x \in L_0}$.\\
Since generic model sets are a broad topic, we restrict ourself
to a more manageable subclass in the next section. It should
also be emphasised that we only consider model sets with
physical space $\Real^2$, for reasons pointed out before.

\subsection{Histogram statistics}\label{sec:histo_stats}

It seems natural to compute statistical data (like variance and skewness) to
analyse the histogram data. We choose not to do so, since this can be misleading.
One can see from the explicit density function $g(t)$ of the $\Integer^2$
case in Sec.~\ref{sec:int_lattice} that the moments of order ${k \ge 2}$
fail to exist. A Taylor expansion gives
\[ g(1/t) = \frac{36}{\pi^4} t^3 + \frac{162}{\pi^6} t^4 + \mathcal{O}(t^5)
   \; \text{ for } \; t \rightarrow 0_{+} \; \text{,} \]
characterising the decay behaviour of the tail. Instead of the statistics,
which just exist because of finite size effects, we provide the coefficients
$c_k$ of $t^k$ (usually two) when the tail of the respective histogram can be
fitted with a power law.

\section{Cyclotomic model sets}\label{sec:cms}

As stated above, we are interested in model set configurations which
admit local visibility tests. This special case is given by
the planar cyclotomic model sets of order ${n \in \Natural}$.
It corresponds to choosing ${d = 2}$, ${k = \phi(n)-2}$ and $L = \CyI{n}$ in
Figure \ref{fig:general_cut_and_project}. Since $\CyI{n} = \CyI{2n}$
for $n$ odd, we impose the condition $n \not\equiv 2$ mod $4$; compare
\cite[Ch.\,3.4]{BaGr13}.\\
The setting can now be used to generate $n$-fold (rotationally) symmetric
point sets (and tilings). The $\star$-map, which maps from physical
to internal space, is given by the extension of an algebraic conjugation;
see \cite{BaGr13} for details.\\
Since the cases ${n = 3, 4}$ yield a planar lattice, we only consider
the configurations with ${n \ge 5}$. Of particular interest are
integers $n$ which admit a simple window test. There are three unique
cases where the window lives in $\Real^2$, or stated differently
where ${\phi(n) = 4}$ holds: 5, 8 and 12.\\

\begin{algorithm2e}
 \SetAlgoLined
 \SetKwInOut{Input}{input}\SetKwInOut{Output}{output}

 \Input{maxsteps, initpoint}
 \Output{vertexlist}

 initialize vertexlist and add initpoint\;
 \For{step $\leftarrow$ $1$ \KwTo maxsteps}{
   \ForEach{$\text{p} \in \text{vertexlist}$}{

     \For{$k \leftarrow 0$ \KwTo ${n-1}$}{
       pp $\leftarrow$ p + $\zeta_n^k$\;

       \uIf{pp is already in vertexlist}{
         skip\;
       }
       \uIf{$\text{pp}^{\star}$ not in window}{
         skip\;
       }

       add pp to vertexlist\;
     }
   }
 }
 \caption{Patch generation for the cyclotomic case.}
 \label{code:patch_generation_cps}
\end{algorithm2e}
The pseudo code in Algorithm \ref{code:patch_generation_cps} then produces the vertices
of a $k$-gon-shaped ($k \in \{10, 8, 12\}$) patch of the corresponding tiling. Note that
for $n=5$, the shape is $10$-fold symmetric because of the $n \not\equiv 2$ (see above)
condition. This $k$-gon shape is desirable because it is already close to being circular and
needs just minor trimming.

\subsection{Ammann--Beenker tiling}\label{sec:amb}

We employ the \emph{Ammann--Beenker} (\textsf{AB}) tiling in its classic version
\cite{AGSh92,BaGr13} with a triangle and a rhombus. It admits a \emph{stone inflation}
(essentially a rule which can be implemented as blowing up the tile followed by a
dissection process), where the triangle (here called the prototile of
\emph{type A}) is inflated as given below in Figure \ref{fig:ab_inflrule_tile_a}:
\putfigure{ab_inflrule_tile_a}
The triangle appears in the tiling with both chiralities, and the other chirality
just uses the reflected rule. The rhombus (prototile of \emph{type B}) appears without
chirality and is inflated according to the rule in Figure \ref{fig:ab_inflrule_tile_b}.
\putfigure{ab_inflrule_tile_b}
Here, the inflation multiplier is given by the \emph{silver mean}
$\lambda_{\text{sm}} = 1 + \sqrt{2}$, which is a \emph{Pisot-Vijayaraghavan}
(\emph{PV}) unit. PV numbers are algebraic integers ${\lambda > 1}$ such that
all algebraic conjugates (except for $\lambda$ itself) lie in the open unit
disk. There is a relation between the regularity of the tiling and the
properties of the inflation multiplier. PV inflations seem to admit more
regular tiling structures \cite[Ch.\,2.5]{BaGr13}; compare
Sec.~\ref{sec:non_pisot} for an example of a less regular tiling point set.\\
A nice property of the \textsf{AB} tiling is that it can be described
as a cyclotomic model set \cite[Ex.\,7.8]{BaGr13}. It corresponds to the diagram in Figure
\ref{fig:general_cut_and_project} of cyclotomic type with parameter
${n = 8}$. The tiling vertices can therefore be described as the set
\[ T_{\textsf{AB}} = \{ x \in \CyI{8} \; : \; x^{\star} \in W_8 \} \text{,} \]
where the $\star$-map is given by the extension of $\zeta_8 \mapsto \zeta_8^3$
and the window $W_8$ is a regular octagon centered at the
origin (edge length one, see Figure \ref{fig:ab_visible_tiling} for the orientation).

The maximal real subring of $\CyI{8}$ is $\Integer[\sqrt{2}]$, with the
unit group generated by ${\pm\lambda_{\text{sm}}}$ from above. By inspecting the action of
these units on the elements of the  $\Integer$-module, one can derive a local visibility test
\[ V_{\textsf{AB}} = \{ x \in T_{\textsf{AB}} \; : \; \lambda_{\text{sm}} x^{\star}
   \notin W_8 \text{ and } x \text{ is coprime} \} \]
for the reference point chosen as the origin. By coprimality of $x$ we
mean coprimality of the coordinates in the direct-sum representation
\[ \CyI{8} \; = \; \Integer[\sqrt{2}] \oplus \Integer[\sqrt{2}] \cdot \zeta_8 \; \text{.} \]
Consider an element $x_1 + x_2 \cdot \zeta_8$ in the above decomposition. The
module $\Integer[\sqrt{2}]$ is a \emph{Euclidean} domain and therefore admits an
algorithm to compute the $\Integer[\sqrt{2}]$-gcd of $x_1$ and $x_2$. By
\emph{coprime} we then understand that this gcd $y$ is a unit, which is equivalent
to $\abs{\algNorm(y)} = 1$, with $\algNorm$ the algebraic norm in the corresponding
module, here given by the map $\algNorm(a + b \cdot \sqrt{2}) = a^2 - 2 \cdot b^2$.
\putfigure{ab_visible_tiling}
The first part of the visibility condition ${x^{\star} \in W_8}$ translates
to the following geometric condition in internal space: If a vertex is visible, then
it lives on a \emph{belt} in internal space, which results from cutting out a
scaled down version of the window from the original window. Both windows are
indicated on the right hand side of Figure \ref{fig:ab_visible_tiling}.
\begin{table}[H]
  \caption{Visibility statistics for the symmetric
           \emph{Ammann--Beenker} tiling.}
  \begin{tabular}{c|c|c|c}\toprule
    maxsteps  & vertices  & visible   & percentage\\
    \midrule
    40        & 561       & 327       & 58.2\%\\
    400       & 47713     & 27561     & 57.7\%\\
    1500      & 662265    & 382221    & 57.7\%\\
    2500      & 1835941   & 1059753   & 57.7\%\\
    \bottomrule
  \end{tabular}
  \label{tab:ab_visibility_stats}
\end{table}
We see that the histogram (generated from roughly $1.8 \cdot 10^6$
vertices) features several characteristics which we have already
observed for the $\Integer^2$-case: A pronounced gap is present
where the distribution has zero mass; then, we have a middle section
where the bulk of the mass is concentrated, and finally a tail
section with a power law decay.
\putfigure{ab_histogram_full}
For an overview of the histogram statistics, see Table
\ref{tab:histogram_stats_cyclotomic} at the end of Sec.~\ref{sec:cms}.

\subsection{T\"ubingen triangle tiling}\label{sec:tueb}

\hspace*{1.2ex}The \emph{T\"ubingen triangle} tiling (\textsf{TT}) is a decagonal case of a
cyclotomic model set with planar window (see \cite{BKSZ90a,BKSZ90b} and
\cite[Ex.\,7.10]{BaGr13}). The underlying module is $\CyI{5}$ with maximal
real subring $\Integer[\tau]$, where $\tau$ is again the multiplier for
the corresponding inflation rule (see Figure \ref{fig:ttt_inflrule_tile_a}
and \ref{fig:ttt_inflrule_tile_b}). See below for a circular patch generated
from applying the inflation rule four times:
\putfigure{ttt_circpatch_infl4}
For the computation of the vertices used for the radial projection, again the
model set description
\[ T_{\textsf{TT}} \; = \; \{ x \in \CyI{5} \; : \; x^{\star} \in W_{10} + \epsilon \} \]
was employed. The window $W_{10}$ is a decagon with edge length $\sqrt{(\tau + 2)/5}$, and
like the \textsf{AB} window, the right-most edge is perpendicular to the $x$-axis.
Here, the $\star$-map is the extension of $\zeta_5 \mapsto \zeta_5^2$.
In this case, we need to apply a small generic shift $\epsilon$ to the window, otherwise
leading to \emph{singular vertices} (vertices which lie on the boundary of the window
when projected to internal space). These are difficult to handle because of precision
issues when testing on the boundary. We therefore restrict ourself to non-singular sets.
In our case we use  $\epsilon = 10^{-4} \cdot (1,1)$ as the shift. The important aspect
here is not to shift in the direction of the window edges. Similar to the eightfold
case, a local visibility condition
\[ V_{\textsf{TT}} \; = \; \{ x \in T_{\textsf{TT}} \; : \; \tau x^{\star} \notin
   W_{10} - \epsilon \text{ and } x \text{ is coprime} \} \]
can be derived. The direct-sum represention here is
$\CyI{5} = \Integer[\tau] \oplus \Integer[\tau] \cdot \zeta_5$, and $\Integer[\tau]$
is again Euclidean.\\
Evaluation with a large patch ($\approx 1.5 \cdot 10^6$ vertices)
produces the following histogram:
\putfigure{ttt_histogram_full}
While being similar to the \textsf{AB} histogram in overall shape, there are numerous
differences in detail, especially in the middle section, which features a lot more
structure and is also nicely aligned to the $\Integer^2$ density function.\\
Zooming into the gap area might even suggest that the middle section decomposes
into smaller components (first step: $(0.18, 0.3)$, second step: $(0.3, 0.5)$, third step:
$(0.5, 1.3)$).
\putfigure{ttt_histogram_zoom}
Again, the statistics can be found in Table \ref{tab:histogram_stats_cyclotomic} below.\\
A related example of a distribution in closed form, for the \emph{golden L} (which is
not a tiling system), has recently been described by Athreya et al. \cite{ACLe13}. It bears
strong resemblence with Fig.~\ref{fig:ttt_histogram_full}, thus making it fall into our
``ordered regime''. This supports the existence of universal features in this approach.

\subsection{G\"ahler's shield tiling}\label{sec:gst}

The \emph{G\"ahler shield} (\textsf{GS}) tiling \cite[Ch.\,5]{Gaeh05} is our last
cyclotomic model set with internal space $\Real^2$. It uses a dodecagonal configuration
\cite[Ex.\,7.12]{BaGr13} and is also interesting in its algebraic properties, which
make the visibility test slightly more involved. The vertex set is
\[ T_{\textsf{GS}} \; = \; \{ x \in \CyI{12} \; : \; x^{\star} \in W_{12} + \epsilon \} \]
with the $\star$-map defined by $\zeta_{12} \mapsto \zeta_{12}^5$. The window $W_{12}$ is
a dodecagon with edge length one and the usual orientation. Again, a shift has to be applied
to avoid singular vertices. The underlying $\Integer$-module decomposes into
\[ \Integer[\sqrt{3}] \oplus \Integer[\sqrt{3}] \cdot \zeta_{12} \; \text{ with } \;
   \lambda_{12} := 2 + \sqrt{3} \]
generating the unit group of $\Integer[\sqrt{3}]$.\\
The local visibility test behaves in a more complex fashion here. Consider an
$x \in \CyI{12}$ and denote by $\algNorm$ the algebraic norm of $\Integer[\sqrt{3}]$.
Now write $x$ in the direct-sum decomposition $x = x_1 + x_2 \cdot \zeta_{12}$ and
define the map
\[ n: \CyI{12} \rightarrow \Natural_1 \; \text{ via } \;
   x \mapsto \abs{N(\gcd(x_1, x_2))} \text{.} \]
Within our finite patch $P$, the set of visible points can then be described as
\begin{align*}
  V_{\textsf{GS}} \; = \; &\{ x \in T_{\textsf{GS}} \; : \; n(x) = 1 \wedge
                            \lambda_1 x^{\star} \notin W_{12} + \epsilon \} \cup {} \\
                          &\{ x \in T_{\textsf{GS}} \; : \; n(x) = 2 \wedge
                            \lambda_2 x^{\star} \notin W_{12} - \epsilon \} \; \text{,}
\end{align*}
where $\lambda_1 := \sqrt{\lambda_{12} \cdot 2}$ and $\lambda_2 := \sqrt{\lambda_{12} / 2}$
(therefore $\lambda_1 \cdot \lambda_2 = \lambda_{12}$), and as long as $\epsilon$ is small
enough in relation to the distances within $P^{\star}$. The first set-component of
$V_{\textsf{GS}}$ is again comprised of \emph{coprime} elements. The second set, however, is
exceptional, and its existence is linked to the degree of the underlying cyclotomic
field, which is ${n = 12}$ here -- a composite number instead of a prime power as in
the other two cases (for cyclotomic fields see \cite{Wash97}).
The difficulty can also be seen on the level of $\Rational(\zeta_n)$, where the unit
group is slightly larger than in the prime power cases, here enlarged by an
additional generating element $z = \sqrt{2 + \sqrt{3}} \cdot \zeta_{24}$.
\putfigure{gst_visible_tiling}
We can see on the right hand side of Figure \ref{fig:gst_visible_tiling} that
two \emph{belts} develop in internal space, one for the coprime vertices and another
one for the exceptional ones. Coprime vertices are represented as
\textcolor{lgray}{\textbf{grey}} dots and exceptional vertices as \textbf{black}
dots. The boundaries of the rescaled (with the factors $\lambda_1$ and
$\lambda_2$ respectively) windows use the same coloring.
\putfigure{gst_histogram_full}
While still retaining the known three-fold structure of the two other cases,
the \textsf{GS} tiling seems to approach the slope-like characteristic
from the Poisson case.
\begin{table}[H]
  \caption{Statistical data generated from the
           radial projection (mean is always $1.0$).}
  \begin{tabular}{c|c|c|c|c|c}\toprule
    tiling       & gap size  & $c_3$  & $c_4$  & $e$   & $k$   \\
    \midrule
    $\Integer^2$ & 0.304     & 0.369  & 0.168  & ---   & ---   \\
    \textsf{AB}  & 0.222     & 0.248  & 0.496  & 2.79  & 38560 \\
    \textsf{TT}  & 0.182     & 0.239  & 0.513  & 2.60  & 31376 \\
    \textsf{GS}  & 0.152     & 0.232  & 0.547  & 4.75  & 67524 \\
    \bottomrule
  \end{tabular}
  \label{tab:histogram_stats_cyclotomic}
\end{table}
The power law fitting was done for the tail starting at $3.0$ (see
Sec.~\ref{sec:histo_stats} for definitions). We indicate the quadratic error by
$e$ in units of $10^{-10}$ and the amount of data points by $k$.

\section{A non-Pisot inflation}\label{sec:non_pisot}

We have seen that the examples of Sec.~\ref{sec:cms} are
qualitatively close to the order properties of the $\Integer^2$
lattice. A similar behaviour of cyclotomic model sets can also
be seen in the mildly related case of discrete tomography \cite{HuSp11}.
One might guess that all kind of deterministic aperiodic tilings
behave that way. However, it turns out that this is not the case.\\
The chiral \emph{Lan\c{c}on--Billard} (\textsf{LB}) tiling \cite{LaBi88} is an
example of an inflation-based tiling with a non-PV multiplier given by
\[ \lambda_{\textsf{LB}} \; = \; \sqrt{\frac{1}{2} \big(5 + \sqrt{5}\big)} \; \text{.} \]
The inflation rule applies to two rhombic prototiles (see
Fig.~\ref{fig:chir_inflrule_tile_a} and Fig.~\ref{fig:chir_inflrule_tile_b}).
\putfigure{chir_inflrule_tile_a}
The resulting tiling vertices live in $\CyI{5}$ (see \cite[Ch.\,6.5.1]{BaGr13}
for details, also concerning the non-PV property of $\lambda_{\textsf{LB}}$), like the \emph{T\"ubingen triangle} tiling above.
\putfigure{chir_inflrule_tile_b}
The \textsf{LB} tiling admits no model set description and it fails to
be a stone inflation, as one can see from the above rules.\\
By multiple inflation of tile A, one can isolate a legal patch
of circular shape that is comprised of five tiles of type A.
We use this patch as our initial seed to grow suitable patches.
\putfigure{chir_circpatch_infl4}
The resulting patches are $C_5$ symmetric and begin to show a high
amount of spatial fluctuation when increasing the number of inflation
steps (the histogram in Figure \ref{fig:chir_histogram_full} was
computed after applying $12$ inflations).
\putfigure{chir_histogram_full}
While not exactly matching the exponential distribution from the
Poisson case, the radial projection at least is sensitive to the higher amount
of spatial disorder in this tiling. In particular, it shows an
exponential rather than a power law decay for large spacings. For the
histogram statistics, see Table \ref{tab:histogram_stats_other} below.

\section{Other planar tilings}

The tilings considered in Secs.~\ref{sec:cms} and \ref{sec:non_pisot}
indicate that the method gives at least partial information about the
order of the point set. Let us look at some more examples.\\
The \emph{chair} tiling \cite{GrSh87} is an example of a inflation tiling
with integer multiplier. It works with just one L-shaped prototile and
can produce patches with $D_4$ symmetry.\\
The patches can also be described as model sets \cite{BaGr13}, but with
a more complicated internal space. We thus employ the inflation method here.
\putfigure{chair_histogram_full}
The vertex set is a subset of $\Integer^2$. It gives a good example why
one has to be careful with the visibility test. Although the set lives
in $\Integer^2$, the standard $\gcd$-test fails in this situation. Consider
a vertex ${p := (x, y)}$ which is not coprime, say with $\gcd(x, y) = k > 1$.
For the integer lattice, one knows that ${p_0 := (\frac{x}{k}, \frac{y}{k})}$
is an element of the set and therefore occludes $p$. This does not need to be
the case here and Figure \ref{fig:badchair_histogram_full} shows that the
difference is indeed significant.
\putfigure{badchair_histogram_full}
The \emph{Penrose--Robinson} (\textsf{PR}) tiling is similar to the \textsf{TT}
on the level of the inflation rule. It uses the same prototiles, but a different
dissection rule \cite[Ch.\,6.2]{BaGr13} after blowing up the tiles by the
inflation factor $\tau$.
\putfigure{prt_histogram_full}
Even though it shares these features with the \textsf{TT}, the resulting
distribution is rather different and offers a high amount of
structure in the bulk section, which can be identified as
\emph{plateau-like} increments.
\putfigure{prt_histogram_zoom}
Another tiling of \emph{Penrose}-type can again be implemented by using a model set
description. This rhombic Penrose (\textsf{RP}) tiling \cite{BKSZ90b} is special in
that it uses a multi-window configuration \cite[Ex.\,7.11]{BaGr13}. Here, the CPS in
Figure \ref{fig:general_cut_and_project} is fixed, but multiple windows $W_i$ are
used. Define the homomorphism
\[ \kappa : \CyI{5} \rightarrow \Integer / 5 \Integer \quad \text{ by } \quad
   \kappa(\sum_{i}{c_i \zeta_5^i}) = \sum_{i}{c_i} \hspace*{-1.0ex}\mod 5, \]
then the window $W_i$ for which the vertex $x \in \CyI{5}$ is tested, is chosen
depending on $\kappa(x)$.
\putfigure{rp_histogram_full}
However, the patches for this case had to be generated using the geometric
visibility test. Although the vertices coming from different
$W_i$ are disjoint, there is still occlusion between the sets which
renders the local test ineffective in this setup.
\begin{table}[H]
  \caption{Statistical data for the other considered
           tilings ($e$ is the error term).}
  \begin{tabular}{c|c|c|c|c|l}\toprule
    tiling       & gap size & $c_3$  & $c_4$   & $c_5$ & \hspace{5ex} $e$   \\
    \midrule
    \textsf{LB}  & 0.0030   & ---    & ---     & ---   & \hspace{4.3ex} --- \\
    chair        & 0.2536   & 0.229  & 0.538   & ---   & $5.07 \cdot 10^{-10}$ \\
    \textsf{PR}  & 0.0783   & 0.066  & 1.339   & ---   & $1.81 \cdot 10^{-10}$ \\
    \textsf{RP}  & 0.1169   & 0.459  & -2.432  & 8.395 & $1.41 \cdot 10^{-9}$  \\
    \bottomrule
  \end{tabular}
  \label{tab:histogram_stats_other}
\end{table}
For the fit of the \textsf{RP} tiling, an additional power was used, to achieve a
similarly small error as in the other cases. Also, a logarithmic fit provides
numerical evidence that the decay behaviour of the chiral \textsf{LB} tiling is
identical to the Poisson case.\\
Another aspect, which is numerically plausible, is the continuous dependence of
the spacing distribution of the cyclotomic cases (Sec.~\ref{sec:cms}) under small
perturbations of the window which leave the area fixed. Replacing the window with
a circle of the same area does not have any noticeable influence on the histogram.
This is in line with related continuity results in \cite{MkSt14} and certainly a
much stronger property than the invariance under removal of singular vertices
(see Sec.~\ref{sec:tueb}), which are known to have density zero in the limit.

\section{Concluding remarks}\label{sec:conclusions}

It would be interesting to study tilings which feature even higher
rotational symmetry than the examples we considered here. While the
data gathered from the three \emph{simple} cyclotomic cases already
shows a tendency, more tilings are needed to fill the picture. The
\emph{de Bruijn} method \cite{deBr81} via dualisation of a grid appears
to be a suitable candidate to generate these kind of tilings.\\
Another aspect which needs further investigation is the existence of a
gap in all studied cases, except the \textsf{LB} one. For cyclotomic
model sets, this seems to be related to the existence of lines with
high density of points on them \cite{Plea03}. This is a feature that
is shared with the $\Integer^2$ case. This has also been observed
in \cite{MkSt14}.\\
Also of interest, but still unclear, is an extension of this
method to higher dimension. A possible way for $\Real^3$ would
be to again project vertices of our set onto the $3$-dimensional
ball of radius $r$. For each projected point $p$, one could now
select the neighbour $q$ with minimal distance to $p$ on the sphere
and consider the angle of the arc between $p$ and $q$. This again produces a list
of angles with which we proceed in the usual way. From a computational point
of view, this case is a lot more involved, since it requires
an exhaustive search for each projected point to find its neighbour.\\
Before closing, we want to point out that projecting from a centre of
maximal symmetry might seem intuitive at first, but still is kind
of special. Since shifting the center indeed changes the distribution,
we want to investigate if some averaging (similar to the shelling
problem \cite{BGJR99} and as also discussed in \cite{MkSt14}) makes
more sense here.

\bibliography{iucr}{}
\bibliographystyle{plain}

\end{document}